\documentclass{amsart}
\usepackage{a4wide}
\usepackage{epsfig,amsmath,amsfonts,amssymb,amsthm,rotating,mathtools}
\usepackage{tikz,pgfplots}
\usetikzlibrary{arrows,snakes,backgrounds}
\usepackage{comment}
\usepackage{algorithm}			
\usepackage{algpseudocode}
\usepackage{mathrsfs}
\usepackage{placeins}
\usepackage{multirow}
\usepackage{hyperref}

\newcommand{\E}{\mathbb{E}}

\newcommand{\isdef}{\mathrel{\mathrel{\mathop:}=}}

\newcommand{\Cor}{\operatorname{Cor}}

\newcommand{\spn}{\operatorname{span}}
\newcommand{\supp}{\operatorname{supp}}

\newcommand{\dist}{\operatorname{dist}}
\newcommand{\diam}{\operatorname{diam}}

\newcommand{\rank}{\operatorname{rank}}
\renewcommand{\div}{\operatorname{div}}

\newcommand{\de}{\operatorname{d}\!}

\newcommand{\balpha}{{\boldsymbol{\alpha}}}
\newcommand{\bbeta}{{\boldsymbol{\beta}}}

\DeclareMathOperator{\Id}{\operatorname{Id}}
\DeclareMathOperator{\tvec}{\operatorname{vec}}

\newcommand{\hsum}{{\text{\tt+}}}
\newcommand{\hdot}{{\text{\tt*}}}

\renewcommand{\d}{\operatorname{d}\!}

\theoremstyle{plain}             
\newtheorem{theorem}{Theorem}[section]

\newtheorem{definition}[theorem]{Definition}

\begin{document}
\title[\(\mathcal{H}\)-matrix based second moment analysis for rough random fields]
{\(\mathcal{H}\)-matrix based second moment analysis for rough random fields
and finite element discretizations}

\author{J.~D\"olz}
\author{H.~Harbrecht}
\author{M.~D.~Peters}
\address{Departement Mathematik und Informatik, Fachbereich Mathematik,
              Universit\"at Basel, \mbox{Spiegelgasse 1}, 4051 Basel, Switzerland}
\email{\{Juergen.Doelz, Helmut.Harbrecht, Michael.Peters\}@unibas.ch}
\thanks{This research has been supported by the
Swiss National Science Foundation (SNSF) through the project ``$\mathcal{H}$-matrix 
based first and second moment analysis''.}

\begin{abstract}
We consider the efficient solution of strongly elliptic partial differential
equations with random load based on the finite element method.
The solution's two-point correlation can efficiently be approximated
by means of an $\mathcal{H}$-matrix, in particular if the correlation 
length is rather short or the correlation kernel is non-smooth. Since the 
inverses of the finite element matrices which correspond to the 
differential operator under consideration can likewise efficiently 
be approximated in the $\mathcal{H}$-matrix format, we can solve 
the correspondent $\mathcal{H}$-matrix equation in essentially linear 
time by using the $\mathcal{H}$-matrix arithmetic. Numerical 
experiments for three-dimensional finite element discretizations 
for several correlation lengths and different smoothness are 
provided. They validate the presented method and demonstrate 
that the computation times do not increase for non-smooth or shortly 
correlated data.
\end{abstract}

\subjclass[2010]{60H35, 65C30, 65N30}

\maketitle

\section{Introduction}
A lot of problems in science and engineering can be modeled in terms of strongly
elliptic boundary value problems. While these problems are numerically well
understood for input data which are given exactly, these input data are often only 
available up to a certain accuracy in practical applications, e.g.,\ due to measurement 
errors or tolerances in manufacturing processes. In recent years, it has therefore 
become more and more important to take these inaccuracies in the input data 
into account and model them as random input parameters.

The Monte Carlo approach, see, e.g.,~\cite{Ca} and the references therein,
provides a straightforward approach to deal with these random data, but
it has a relatively slow convergence rate which is only in the sense of the root
mean square error. This, in turn, means that a large amount of samples has to 
be generated to obtain computational results with an acceptable accuracy, whereas 
the results still have a small probability of being too far away from the true solution.
Therefore, in the past several years there have been presented multiple
deterministic approaches to overcome this obstacle. For instance, random
loads have been  considered in \cite{ST03a,vPS06}, random coefficients
in \cite{BNT07,BTZ04,DBO,FST05,GS91,KS11,MK05,NTW08}, and random 
domains in \cite{HSS08a,XT06}.

For a domain \(D\subset\mathbb{R}^d\) and a
probability space \((\Omega,\mathcal{F}, \mathbb{P})\),
we consider the Dirichlet problem
\begin{equation}\label{eq:modprob}
\left. \begin{aligned}
\mathcal{L}u(\omega,{\bf x}) &= f(\omega,{\bf x})&&\text{ for }{\bf x}\in D\\
u(\omega,{\bf x}) &= 0 &&\text{ for }{\bf x}\in\Gamma\isdef\partial D\quad
\end{aligned}\right\}\quad\text{\(\mathbb{P}\)-almost surely}
\end{equation}
with random load $f(\omega,{\bf x})$ and a strongly elliptic partial
differential operator $\mathcal{L}$ of second order.

We can compute the solution's mean
\[
\mathbb{E}_u({\bf x})\isdef\int_\Omega u(\omega,{\bf x})\de\mathbb{P}(\omega)
\]
and also its two-point correlation
\[
\Cor_u({\bf x},{\bf y})\isdef\int_\Omega u(\omega,{\bf x})u(\omega,{\bf y})
\de\mathbb{P}(\omega)
\]
if the respective quantities of the input data are known. 
Namely, the mean \(\E_u\) satisfies
\begin{equation}\label{eq:meanpdei}
\mathcal{L}\E_u=\E_f\text{ in }D\quad\text{and}\quad \E_u=0\text{ on }\Gamma
\end{equation}
due to the linearity of the expectation and the differential operator 
\(\mathcal{L}\). Taking into account the multi-linearity 
of the tensor product, one verifies by tensorizing 
\eqref{eq:modprob} that
\begin{equation}\label{eq:coreq}
\begin{aligned}
(\mathcal{L}\otimes\mathcal{L})\Cor_u&=\Cor_f&&\text{ in }D\times D,\\
(\mathcal{L}\otimes\Id) \Cor_u &= 0&&\text{ on }D\times\Gamma,\\
(\Id\otimes\mathcal{L}) \Cor_u &= 0&&\text{ on }\Gamma\times D,\\
\Cor_u&=0&&\text{ on }\Gamma\times\Gamma.
\end{aligned}
\end{equation}
From $\Cor _u$, we can compute the variance $\mathbb{V}_u$ 
of the solution due to
\[
\mathbb{V}_u (\mathbf{x}) = \Cor _u(\mathbf{x},\mathbf{x})-\mathbb{E}_u(\mathbf{x})^2.
\]

If a low-rank factorization of $\Cor _f$ is available,
\eqref{eq:coreq} can easily be solved by standard finite element
methods. The existence of an accurate low-rank approximation is
directly related to the spectral decomposition of the associated 
Hilbert-Schmidt operator
\begin{equation}\label{eq:HS}
  (\mathcal{C}_f \psi)({\bf x})\isdef\int_D \Cor_f({\bf x},{\bf y})\psi({\bf y})\d{\bf y}.
\end{equation}
Let $\Cor _f\in H^p(D)\otimes H^p(D)$, then, according 
to \cite{GH14,ST06}, the eigenvalues of this Hilbert-Schmidt 
operator decay like 
\begin{equation}\label{eq:eigdec}
\lambda _{m}\lesssim m^{-2p/d}\text{ as }m\rightarrow\infty.
\end{equation}
Unfortunately, the constant in this estimate behaves similar to the $H^p(D)\otimes 
L^2(D)$-norm of $\Cor _f$. The following consideration shows that
this can lead to large constants in the decay estimate if the correlation length
is small. Let the correlation kernel $k(r)$ depend only on the distance 
$r = \|{\bf x}-{\bf y}\|$. Then, the derivatives 
$\partial _{\mathbf{x}}^{\boldsymbol\alpha}\Cor _f(\mathbf{x},\mathbf{y})$
and
$\partial _{\mathbf{y}}^{\boldsymbol\alpha}\Cor _f(\mathbf{x},\mathbf{y})$
of the correlation
\[
\Cor _f(\mathbf{x},\mathbf{y})=k\bigg(\frac{\|\mathbf{x-\mathbf{y}\|}}{\ell}\bigg),
\]
involve the factor $\ell ^{-|\boldsymbol\alpha|}$, leading to a 
constant $\ell ^{-p}$ in the decay estimate of the eigenvalues
\eqref{eq:eigdec}. Thus, for a small correlation length $\ell$, a 
low-rank approximation of $\Cor _f$ becomes prohibitively 
expensive to compute.

Different approaches to tackle the solution of \eqref{eq:coreq} 
have been considered in several articles, where most of them have
in common that they are in some sense based on a sparse tensor 
product discretization of the solution. For example, the computation 
of the second moment, i.e.,~\(\Cor_u\), has been considered for elliptic 
diffusion problems with random loads in \cite{ST03a} by means of a 
sparse tensor product finite element method. A sparse tensor product 
wavelet boundary element method has been used in \cite{HSS08a} 
to compute the solution's second moment for elliptic potential 
problems on random domains. In \cite{Har10,HSS08b}, the 
computation of the second moment was done by multilevel finite 
element frames. Recently, this concept has been simplified by 
using the combination technique, cf.~\cite{HPS13a}. Unfortunately,
the sparse tensor product discretization needs to resolve the concentrated 
measure for short correlation lengths. This means that the number of
hierarchies of the involved finite element spaces has to be doubled
if the correlation length is halved to get the same accuracy. 

The present article discusses a different approach to approximate 
the full tensor product discretization without losing its resolution 
properties. In \cite{DHP15}, it has been demonstrated that the 
$\mathcal{H}$-matrix technique is a powerful tool to cope with Dirichlet 
data of low Sobolev smoothness if the problem is formulated as a 
boundary value problem. There, the similar behavior of two-point 
correlation kernels and boundary integral operators has been exploited.
In \cite{DHS15}, $\mathcal{H}$-matrix compressibility of the solution 
was proven also in case of local operators on domains. In the present 
article, we will combine this theoretical foundation with the 
$\mathcal{H}$-matrix technique used in \cite{DHP15} to efficiently 
solve \eqref{eq:coreq} by the finite element method for a right hand side
$\Cor _f$ with small correlation length or low Sobolev smoothness.

The general concept of \(\mathcal{H}\)-matrices and the corresponding 
arithmetic have at first been introduced in \cite{Hac99,HK00}.
\(\mathcal{H}\)-matrices are feasible for the data-sparse representation 
of (block-) matrices which can be approximated block-wise with low-rank
and have originally been employed for the efficient treatment of boundary integral
equations as they arise in the boundary element method.
 
The rest of this article is organized as follows. In
Section~\ref{sec:preliminaries}, we introduce the Galerkin discretization of
the problem under consideration. Section~\ref{sec:Hmatrixcorrelationkernels} 
discusses the compressibility of discretized correlation kernels and the 
efficient solution of general correlation equations. In Section~\ref{sec:FEMHmatrices},
we recall some specialities of $\mathcal{H}$-matrices in the context of finite elements. 
In particular, we restate a phenomenon, called ``weak admissibility'', which produces
a more data-sparse representation of the correlation matrices, and 
$\mathcal{H}$-matrix nested dissection techniques. In 
Section~\ref{sec:numres}, we present numerical examples to validate 
and quantify the proposed method. Finally, in Section~\ref{sec:conclusion}, 
we draw our conclusions from the theoretical findings and the numerical results.

\section{Preliminaries}\label{sec:preliminaries}
For the remainder of this article, let \(D\subset\mathbb{R}^d\) be a 
Lipschitz domain, \((\Omega,\mathcal{F},\mathbb{P})\) 
a separable, complete probability space and $\mathcal{L}$ the linear 
differential operator of second order given by
\begin{equation}\label{eq:diffop}
(\mathcal{L}u)({\bf x})\isdef -\div\big(\mathbf{A}(\mathbf{x})\cdot\nabla u(\mathbf{x})\big)
+c(\mathbf{x})u(\mathbf{x}).
\end{equation}
The differential operator shall be strongly elliptic
in the sense that
\[
\underline{\alpha}\|\boldsymbol\xi\|_2^2\leq
\boldsymbol\xi^{\intercal}\mathbf{A}(\mathbf{x})\boldsymbol\xi\leq
\overline{\alpha}\|\boldsymbol\xi\|_2^2\text{ for all }\boldsymbol\xi\in\mathbb{R}^d
\]
with coefficients $\mathbf{A}\in W^{1,\infty}\big(D,\mathbb{R}^{d\times d}\big)$,
$0\leq c\in L^{\infty}(D,\mathbb{R})$,
and $0<\underline{\alpha}\leq\overline{\alpha}<\infty$.

Under these assumptions, for a given load $f\in L^{2}_{\mathbb{P}}\big(\Omega ,H^{-1}(D)\big)$,
the Dirichlet problem
\[
\begin{aligned}
\mathcal{L}u(\omega,{\bf x}) &= f(\omega,{\bf x})&&\text{ for }{\bf x}\in D\\
u(\omega,{\bf x}) &= 0 &&\text{ for }{\bf x}\in\Gamma\isdef\partial D,
\end{aligned}
\]
is known to have a unique solution $u(\omega ,\cdot )\in H_0^1(D)$ for
$\mathbb{P}$-almost every $\omega\in\Omega$, cf.,~e.g.,~\cite{GT77}. As 
a result, the mean $\mathbb{E}_u\in H_0^1(D)$ and the correlation
$\Cor _u\in H_0^1(D)\otimes H_0^1(D)$ are well defined.

For the efficient numerical solution of \eqref{eq:coreq},
we use a finite element Galerkin scheme.
To that end, we introduce a finite element space $V_N 
= \spn\{\varphi_1,\ldots,\varphi_N\}\subset H_0^1(D)$. It is 
assumed that the mesh which underlies this finite 
element space is quasi-uniform. 
The basis functions $\{\varphi_i\}_i$ are 
supposed to be locally and isotropically
supported such that 
$\diam(\supp\varphi_i)\sim N^{-1/d}$. In particular, 
we can assign to each degree of freedom $i\in\{1,\ldots,N\}$ 
a suitable point ${\bf x}_i\in D$, e.g.,\
the barycenter of the support of the corresponding 
basis function or the corresponding Lagrangian interpolation
point if nodal finite element shape functions are considered. 

The variational formulation of \eqref{eq:coreq} is given as follows:
\[
\begin{aligned}
  &\text{Find}\ \Cor_u\in H_0^1(D)\otimes H_0^1(D)\text{ such that } \\
  &\quad\big((\mathcal{L}\otimes\mathcal{L})\Cor_u,v\big)_{L^2(D\times D)}
	=(\Cor_f,v)_{L^2(D\times D)}
		\text{ for all }v\in H_0^1(D)\otimes H_0^1(D).
\end{aligned}
\]
By replacing the energy space \(H_0^1(D)\otimes H_0^1(D)\) 
in this variational formulation by the finite dimensional 
ansatz space \(V_N\otimes V_N\), we arrive at
\begin{equation}\label{eq:galerkin}
\begin{aligned}
  &\text{Find}\ \Cor_{u,N}\in V_N\otimes V_N\text{ such that } \\
  &\quad\big((\mathcal{L}\otimes\mathcal{L})\Cor_{u,N},v\big)_{L^2(D\times D)}
	=(\Cor_f,v)_{L^2(D\times D)}
		\text{ for all }v\in V_N\otimes V_N.
\end{aligned}
\end{equation}
A basis in \(V_N\otimes V_N\) is formed by the set of tensor product 
basis functions $\{\varphi_i \otimes\varphi_j\}_{i,j}$.
Hence, 
representing $\Cor_{u,N}$ by its basis expansion, yields
\[
\Cor_{u,N}=\sum_{\ell,\ell'=1}^{N}u_{\ell,\ell'}(\varphi_\ell\otimes\varphi_{\ell'}).
\]
Setting ${\bf C}_u\isdef[u_{\ell,\ell'}]_{\ell,\ell'}$,
we end up with the linear system of equations
\begin{equation}\label{eq:linsys}
({\bf A}\otimes{\bf A})\tvec({\bf C}_u)=\tvec({\bf C}_f),
\end{equation}
where ${\bf C}_f\isdef\big[(\Cor_f,\varphi_\ell\otimes\varphi_{\ell'}
)_{L^2(D\times D)}\big]_{\ell,\ell'}$ is the discretized 
two-point correlation of the Dirichlet data $f$ and ${\bf A}\isdef
\big[(\mathcal{L}\varphi_{\ell'},\varphi_{\ell})_{L^2(D)}\big]_{\ell,\ell'}$ 
is the system matrix of the second order differential operator \eqref{eq:diffop}. 
In \eqref{eq:linsys}, the tensor product has, as usual in connection with
matrices, to be understood as the Kronecker product. Furthermore, for a 
matrix \({\bf B}=[{\bf b}_1,\ldots,{\bf b}_n]\in\mathbb{R}^{m\times n}\), 
the operation $\tvec({\bf B})$ is defined as
\[
\tvec([{\bf b}_1,\ldots,{\bf b}_n])\isdef
\begin{bmatrix}{\bf b}_1\\ \vdots\\ {\bf b}_n\end{bmatrix}\in\mathbb{R}^{mn}.
\]

For matrices \({\bf B}\in\mathbb{R}^{k\times n}\), \({\bf C}\in
\mathbb{R}^{\ell\times m}\) and \({\bf X}\in\mathbb{R}^{m\times n}\),
there holds the relation
\[
({\bf B}\otimes{\bf C})\tvec({\bf X})=\tvec({\bf CXB}^\intercal).
\] 
Hence, we may rewrite \eqref{eq:linsys} according to
\begin{equation}\label{eq:linsys2}
{\bf A}{\bf C}_u{\bf A}^\intercal={\bf C}_f.
\end{equation}

An approach to deal with non-homogeneous boundary conditions 
has been presented in \cite{Har10}.
\section{$\mathcal{H}$-matrix approximation of correlation kernels}
\label{sec:Hmatrixcorrelationkernels}
The matrix equation \eqref{eq:linsys2} has $N^2$ unknowns 
and is therefore not directly solvable if $N$ is large due to memory and 
time consumption. Thus, an efficient compression scheme and a powerful 
arithmetic are needed to obtain its solution. In the following, we will restrict 
ourselves to \emph{asymptotically smooth} correlation kernels $\Cor _f$,
i.e.,~correlation kernels satisfying the following definition.
 
\begin{definition} Let \(k\colon\mathbb{R}^d\times\mathbb{R}^d
\to\mathbb{R}\). The function \(k\) is called \emph{asymptotically 
smooth} if for some constants \(c_1,c_2>0\) and \(q\in\mathbb{R}\) holds
\begin{equation}	\label{eq:assymptoticallysmooth}
  \big|\partial_{\bf x}^{\boldsymbol\alpha}\partial_{\bf y}^{\boldsymbol\beta} 
  	k({\bf x},{\bf y})\big| \leq c_1\frac{(|\boldsymbol{\alpha}|+|\boldsymbol{\beta}|)!}
	{c_2^{{|\boldsymbol\alpha|+|\boldsymbol\beta|}}}
		\|{\bf x}-{\bf y}\|_2^{-d-2q-|\boldsymbol\alpha|-|\boldsymbol\beta|},\quad\mathbf{x}\neq\mathbf{y},
\end{equation}
independently of \(\balpha\) and \(\bbeta\). 
\end{definition}

Examples for asymptotically smooth correlation kernels are the
Mat\'ern kernels, which include especially the Gaussian kernel,
cf.~\cite{Mat86,RW05} and the references therein. A main feature 
of such asymptotically smooth correlation kernels is that they exhibit 
a data-sparse representation by means of $\mathcal{H}$-matrices,
cf.,~e.g.,~\cite{Beb08,Bor10,Hac15}. 

$\mathcal{H}$-matrices rely on local low-rank approximations 
of a given matrix ${\bf X}\in\mathbb{R}^{N\times N}$. For suitable 
non-empty index sets $\nu,\nu'\subset\{1,\ldots,N\}$, a matrix block 
${\bf X}|_{\nu\times\nu'}$ can be approximated by a rank-$k$ matrix. 
This approximation can be represented 
in factorized form ${\bf X}|_{\nu\times\nu'}\approx{\bf YZ}^\intercal$ 
with factors ${\bf Y}\in\mathbb{R}^{\nu\times k}$ and ${\bf Z}\in
\mathbb{R}^{\nu'\times k}$. Hence, if $k\ll\min\{\#\nu,\#\nu'\}$, 
the complexity for storing the block is considerably reduced. 
The construction of the index sets is based on the \emph{cluster tree}.
%
\subsection{Cluster tree}
\label{sec:ClTr}
For a tree $\mathcal{T}=(V,E)$ with vertices $V$ and edges $E$,
we define its set of leaves by
\[
\mathcal{L}(\mathcal{T})\isdef\{\sigma\in V\colon\sigma~\text{has no sons}\}.
\]
Furthermore, we say that $\mathcal{T}$ is a \emph{cluster tree} for
the set $\{1,\ldots,N\}$ if the following conditions hold.
\begin{itemize}
\item $\{1,\ldots,N\}$ is the \emph{root} of $\mathcal{T}$.
\item All $\sigma\in V\setminus\mathcal{L}(\mathcal{T})$
      are the disjoint union of their sons.
\end{itemize}
The \emph{level} of $\sigma\in\mathcal{T}$ is its distance of the root,
i.e.,~the number of son relations that are required for traveling from
$\{1,\ldots,N\}$ to $\sigma$. We define the set of clusters
on level $j$ as
\[
\mathcal{T}^{(j)}:=\{\sigma\in\mathcal{T}\colon \sigma~\text{has level}~j\}.
\]
The construction of the cluster tree is based on the
{\em support} of the clusters. The support $\Upsilon_{\sigma}$ 
of a cluster $\sigma$ is defined as the union of the supports 
of the basis functions corresponding to their elements, that is
 \[
   \Upsilon_{\sigma} = \bigcup_{i\in\sigma}\Upsilon_i\text{ where}\ 
   	\Upsilon_i:=\supp\varphi_i\ \text{for all}\ i\in\{1,\ldots,N\}.
 \]
For computing complexity bounds, the cluster tree
should match the following additional requirements,
uniformly as $N\to\infty$:
\begin{itemize}
\item 
The cluster tree is 
a balanced tree in the sense that the maximal level satisfies 
$J\sim\log_2 N$.
\item 
The diameter of the support $\Upsilon_{\sigma_{j}}$, $\sigma_{j}\in\mathcal{T}^{(j)}$,
is local with respect to the level $j$, i.e.,
$\diam\Upsilon_{\sigma_{j}}\sim 2^{-j/d}$. 
Moreover, 
the number $\#\sigma_{j}$ of indices contained in a cluster 
$\sigma_{j}\in\mathcal{T}^{(j)}$ scales approximately 
like $2^{J-j}$, i.e.,\ $\#\sigma_{j}\sim 2^{J-j}$.
\end{itemize}

Until further notice, a binary
cluster tree $\mathcal{T}$ with the indicated terms 
should be given for our further considerations. A common 
algorithm for its construction is based on a hierarchical 
subdivision of the point set which is associated with the 
basis functions, cf.,~e.g.,~\cite{Beb08,Bor10,Hac15}. 
We begin by embedding the point set 
$\{\mathbf{x}_1,\ldots,\mathbf{x}_N\}$ in a top-level bounding-box. This bounding-box
is subsequently subdivided into two cuboids of the same size
where the corresponding clusters are described by the
points in each bounding-box. This process 
is iterated until a bounding-box encloses less than a 
predetermined number of points. 

\subsection{$\mathcal{H}$-Matrix approximation}
\label{sec:H-matrices}
$\mathcal{H}$-matrices have originally been invented in
\cite{Hac99,HK00} and are a generalization of 
cluster techniques for the rapid solution of boundary 
integral equations such as the fast multipole method 
\cite{GR87}, the mosaic skeleton approximation \cite{Tyr96}, 
or the adaptive cross approximation \cite{Beb00}.

For the discretization of an asymptotically smooth correlation, 
we introduce a partition of its domain of definition which 
separates smooth and non-smooth areas of the kernel function.
It is based on the following
\begin{definition} \label{def:adm}
Two clusters $\sigma$ and $\sigma'$ are called \emph{$\eta$-admissible}
if
\begin{equation}\label{eq:admissibility}
\max\{\diam(\Upsilon_{\sigma}),\diam(\Upsilon_{\sigma'})\}
\le\eta\,\dist(\Upsilon_{\sigma},\Upsilon_{\sigma'})
\end{equation}
holds for some fixed \(\eta>0\).
\end{definition}

We can obtain the set of admissible blocks by means
of a recursive algorithm: Starting with 
the root $(\sigma_{0,0},\sigma_{0,0})$, the bounding-boxes of the current cluster pair are checked 
for admissibility. If they are admissible, the cluster pair is added to the set $\mathcal{F}$ which 
corresponds to the correlation kernel's \emph{farfield}. Otherwise, the 
admissibility check will be performed on all bounding-boxes of the possible pairs of son clusters of
the two original clusters. When we arrive at a pair of leaf
clusters with inadmissible bounding-boxes, the clusters are added to
the set $\mathcal{N}$ which corresponds to the correlation kernel's
\emph{nearfield}. The set
$\mathcal{B}=\mathcal{F} \cup\mathcal{N}$ obviously inherits a tree structure from the recursive construction of $\mathcal{F}$ and
$\mathcal{N}$ and is called the \emph{block cluster tree},
see~\cite{Beb08,Bor10,Hac15}.

With the definition of the block cluster tree at hand, we are 
finally in the position to introduce \(\mathcal{H}\)-matrices.

\begin{definition}
The set $\mathcal{H}(\mathcal{B},k)$ of \emph{$\mathcal{H}$-matrices}
of maximal block rank $k$ is defined according to
\[ 
  \mathcal{H}(\mathcal{B},k):=\big\{{\bf X}\in\mathbb{R}^{N\times N}:
  \rank\big({\bf X}|_{\sigma\times\sigma'}\big)\le k\text{ for all }\big(\sigma,\sigma'\big)\in\mathcal{F}\big\}.
\]
Note that all nearfield blocks ${\bf X}|_{\sigma\times\sigma'}$,
$(\sigma,\sigma')\in\mathcal{N}$, are allowed to be full matrices.
\end{definition}

In accordance with \cite{Beb08,Bor10,Hac15}, the storage cost of 
an $\mathcal{H}$-matrix ${\bf X}\in\mathcal{H}(\mathcal{B},k)$ is 
$\mathcal{O}(k N\log N)$. Here, for asymptotically smooth correlation 
kernels, the rank $k$ depends poly-logarithmically on the desired 
approximation accuracy $\varepsilon$, which in turn depends usually 
on the degrees of freedom $N$. These remarks pertain to the
approximation of an \emph{explicitly} given, asymptotically smooth 
correlation kernel $k$, such as $\Cor_f$ in \eqref{eq:coreq}.

The compressibility of an \emph{implicitly} given correlation kernel,
such as $\Cor_u$ in \eqref{eq:coreq}, has been studied in \cite{DHS15} 
for the case of smooth domains $D$. We restate the main theorem for
the setting of the present article which employs that the Hilbert-Schmidt 
operator \eqref{eq:HS}, related with the correlation kernel $\Cor_f$, 
is in general a \emph{pseudo-differential operator}, see, e.g.,\ 
\cite{Hor03,Hor07,HW08,Tay81} and the references therein.

\begin{theorem}
In the domain $D$ with analytic boundary $\partial D$, 
assume that the correlation kernel $\Cor_f$ in \eqref{eq:coreq}
gives rise to an operator $\mathcal{C}_f \in OPS^\theta_{cl,s}$, 
i.e.,\ to a classical pseudo-differential operator with symbol 
$a_{f}(\mathbf{x},\boldsymbol\xi)$ of order $\theta$ and of Gevrey class $s\geq 1$ 
in the sense of \cite[Def.\ 1.1]{KM67}. Assume further that the 
coefficients of the differential operator $\mathcal{L}$ are smooth.
Then, the correlation kernel $\Cor _u$ of \eqref{eq:coreq}  
is the Schwartz kernel of an operator 
$\mathcal{C} _u\in OPS^{\theta-4}_{cl,s}$.

Moreover, the kernel $\Cor_u(\mathbf{x},\mathbf{y})$ of the correlation operator 
$\mathcal{C}_u$ is smooth in $D\times D$ outside of the diagonal 
$\Delta := \{(\mathbf{x},\mathbf{y})\in D\times D:\mathbf{x}=\mathbf{y}\}$ and there holds the 
pointwise estimate
\begin{equation}
\label{eq:pseudoest}
\qquad\qquad|\partial^{\boldsymbol\alpha}_\mathbf{x} \partial^{\boldsymbol\beta}_\mathbf{y} \Cor_u({\bf x},{\bf y})|
\leq c \mathscr{A}^{|\boldsymbol\alpha+\boldsymbol\beta|} (|\boldsymbol\alpha|!)^s \boldsymbol\beta! \|{\bf x}-{\bf y}\|_2^{-\theta -d-|\boldsymbol\alpha|-|\boldsymbol\beta|+4}
\end{equation}
for all \(\boldsymbol\alpha,\boldsymbol\beta\in \mathbb{N}_0^d\), \(({\bf x},{\bf y})\in(D\times D)\backslash\Delta\) with
some constants $c$ and $\mathscr{A}$ which depend 
only on $D$ and on $a_{f}$.
\end{theorem}

Obviously, for $s=1$, estimate \eqref{eq:pseudoest} directly implies
condition \eqref{eq:assymptoticallysmooth} for the asymptotic smoothness
of $\Cor _u$, allowing us to approximate $\Cor _u$ by the means of
$\mathcal{H}$-matrices. In particular, \cite{DHS15} provides also some
numerical evidence that this result could likely be extended to
Lipschitz domains.

An example of correlation kernels for $\Cor _f$ satisfying the
condition of this theorem for $s\geq 1$ is the Mat\'ern class
of kernels. We refer to \cite{DHS15} for more details on how to
verify the assumptions of the theorem for other correlation kernels.
%
\subsection{$\mathcal{H}$-Matrix arithmetic and iterative solution}
An important feature of $\mathcal{H}$-matrices is that 
efficient algorithms for approximate matrix arithmetic 
operations are available. For two $\mathcal{H}$-matrices ${\bf H}_1,{\bf H}_2
\in\mathcal{H}(\mathcal{B},k)$, the approximate \emph{matrix-matrix 
addition} ${\bf H}_1\hsum {\bf H}_2\in\mathcal{H}(\mathcal{B},k)$
can be performed in $\mathcal{O}(k^2 N\log N)$ operations while the approximate 
\emph{matrix-matrix multiplication} \mbox{${\bf H}_1\hdot {\bf H}_2\in\mathcal{H}(\mathcal{B},k)$} can be performed in $\mathcal{O}
(k^2 N\log^2 N)$ operations. Both of these operations are essentially
block matrix algorithms with successive recompression schemes.
Moreover, employing the recursive 
block structure, the approximate inversion or the approximate 
computation of the $LU$-decomposition within $\mathcal{H}
(\mathcal{B},k)$ can also be performed in only $\mathcal{O}
(k^2 N\log^2 N)$ operations. We refer the reader to 
\cite{Beb08,Bor10,GH03,Hac15,HK00} for further results and 
implementation details.
Especially, the parallelization of the $\mathcal{H}$-matrix arithmetic and
the $\mathcal{H}$-LU-decomposition has been discussed in \cite{Kri05,Kri13}.

In the context of correlation equations, this approximate
$\mathcal{H}$-matrix arithmetic has successfully been used in \cite{DHP15}
to solve a problem, similar to \eqref{eq:coreq}, which has been 
discretized by the boundary element method. Then, the
matrix $\mathbf{A}$ in \eqref{eq:linsys2} corresponds to the stiffness matrix 
from the boundary element method which can naturally be approximated 
by the $\mathcal{H}$-matrix technique. The resulting matrix equation
has been solved using an iterative solver based  on \emph{iterative 
refinement}, cf.~\cite{GVL12,Mol67,Wil63}, which we are also
going to employ here. This method has originally been introduced in 
\cite{Wil63} for the improvement of solutions to linear systems of 
equations based on the LU-factorization. 

Having all matrices in \eqref{eq:linsys2} represented by 
$\mathcal{H}$-matrices, the solution can then be approximated 
as follows. Let ${\bf A}\approx\hat{\bf L}\hat{\bf U}$, where 
$\hat{\bf L},\hat{\bf U}\in\mathcal{H}(\mathcal{B},k)$, be an 
approximate LU-decomposition to \({\bf A}\), e.g.,~computed from 
$\mathbf{A}$ by the $\mathcal{H}$-matrix arithmetic. Starting 
with the initial guess
$
{\bf C}_u^{(0)}=\hat{\bf U}^{-1}\hat{\bf L}^{-1}{\bf C}_f\hat{\bf L}^{-\intercal}\hat{\bf U}^{-\intercal},
$
we iterate
\begin{equation}
\label{eq:iterativerefinement}
{\boldsymbol{\Theta}}^{(i)}= {\bf C}_f-{\bf A}{\bf C}_u^{(i)}{\bf A}^\intercal,
\quad{\bf C}_u^{(i+1)}={\bf C}_u^{(i)}
+\hat{\bf U}^{-1}\hat{\bf L}^{-1}{\boldsymbol{\Theta}}^{(i)}\hat{\bf L}^{-\intercal}\hat{\bf U}^{-\intercal},
\quad i=0,1,\ldots.
\end{equation}
Note that we use, in contrast to \cite{DHP15}, the LU-decomposition 
with forward and backward substitution algorithms which avoids the 
expensive computation of an approximate inverse. Whether we 
use an approximate inverse as in \cite{DHP15} or an approximate
LU-decomposition,
the idea of the iterative refinement stays the same: The residual 
\({\boldsymbol{\Theta}}^{(i)}\) is computed with a higher precision 
than the correction \(\hat{\bf U}^{-1}\hat{\bf L}^{-1}{\boldsymbol{\Theta}}^{(i)}\hat{\bf L}^{-\intercal}\hat{\bf U}^{-\intercal}\). This yields an improved approximation to 
the solution in each step. Note that this algorithm also algebraically 
coincides with an undamped preconditioned Richardson iteration, 
see, e.g.,~\cite{saa03}.

In the following, we will elaborate how this approach can be 
realized in the context of the finite element method.
If $\mathbf{A}$ is symmetric and positive definite, the LU-decomposition
could also be replaced by a Cholesky decomposition. Nonetheless, we will see
in the numerical experiments that the computation time of the decomposition
is negligible compared to the overall computation time and we prefer to stay
in the more general, i.e.~non-symmetric, setting.

\section{$\mathcal{H}$-matrices in the context of finite elements}
\label{sec:FEMHmatrices}
Although a finite element matrix has a sparse structure,
its inverse and both factors of its LU-decomposition
are generally fully populated. Nevertheless, the inverse 
and the LU-decomposition exhibit a data-sparse structure in the sense that 
they are $\mathcal{H}$-matrix compressible. We recall the main 
concepts from the literature, see,
e.g.,~\cite{Beb05,BH03,Fau15,FMP15,Hac15}.

\subsection{General concepts}
A rough argument for the $\mathcal{H}$-matrix compressibility of the inverse
makes use of the Green's function $G$ of $\mathcal{L}$.
Let $\delta_\mathbf{x}$ denote the Dirac distribution at
the point $\mathbf{x}$ and let
$G\colon\mathbb{R}^d\times\mathbb{R}^d\rightarrow\mathbb{R}$
satisfy
\[
\mathcal{L}_{\mathbf{y}}G(\mathbf{x},\mathbf{y})=\delta _{\mathbf{x}}
\text{ and }G(\cdot ,\mathbf{y})|_{\Gamma}=0.
\]
Then, the solution of
\[
\begin{aligned}
\mathcal{L}u({\bf x}) &= f({\bf x})&&\text{ for }{\bf x}\in D,\\
u({\bf x}) &=0 &&\text{ for }{\bf x}\in\Gamma,
\end{aligned}
\]
can be represented by
\[
u(\mathbf{x}) = (\mathcal{L}^{-1}f)(\mathbf{x}) = \int _DG(\mathbf{x},
\mathbf{y})f(\mathbf{y})\de\mathbf{y},\quad\mathbf{x}\in D.
\]
If the Green's function is analytic away from the diagonal, e.g.,~in
the case of constant coefficients of $\mathcal{L}$, we can approximate the Green's
function away from the diagonal by local expansions of the kind
\[
G(\mathbf{x},\mathbf{y})\approx\sum _{i=1}^ka(\mathbf{x})b(\mathbf{y}),
\]
which is the theoretical basis for an $\mathcal{H}$-matrix approximation,
see \cite{Beb08,Bor10,Hac15}.

However, one of the advantages of the finite element method
is that it can be applied also in case of non-constant coefficients. In 
\cite{BH03}, a proof was presented to guarantee the existence
of an $\mathcal{H}$-matrix approximation to the inverse of the finite
element stiffness matrix even in the case of essentially bounded diffusion
coefficients and the other coefficients set to zero. This result was
then extended in \cite{Beb05} to allow all coefficients to be only essentially bounded,
providing the theoretical foundation for an $\mathcal{H}$-matrix approximation
to the inverse of the differential operator from \eqref{eq:diffop}.
Having the $\mathcal{H}$-matrix approximability of the inverse
to the finite element matrix available, the approximability of the
LU-decomposition to the finite element matrix has then been proven in
\cite{Beb07}.

While these first results hold up to the finite element discretization
error, the results have recently been improved in \cite{Fau15,FMP15} to hold
without additional error.

\begin{figure}
\centering
\begin{minipage}{0.32\textwidth}
\centering
\includegraphics[width=0.95\textwidth]{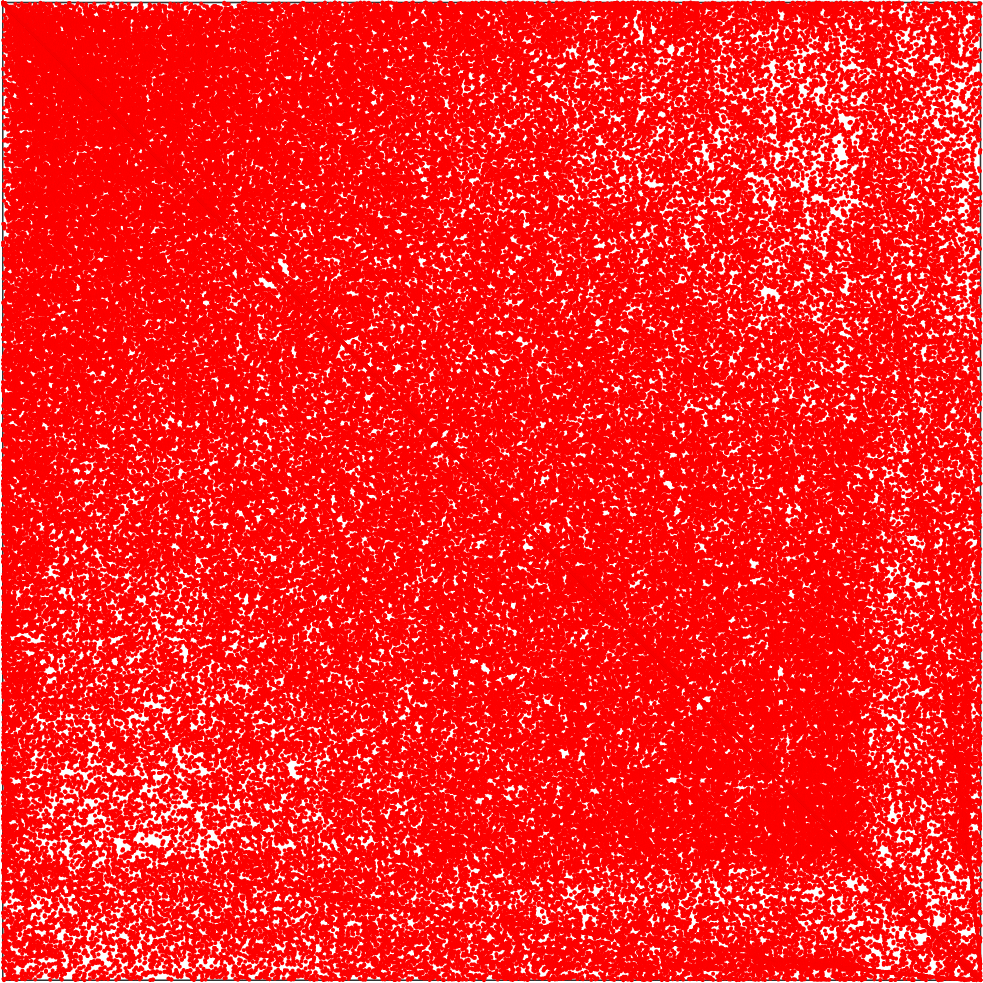}\\
sparse FEM-matrix
\end{minipage}
\begin{minipage}{0.32\textwidth}
\centering
\includegraphics[width=0.95\textwidth]{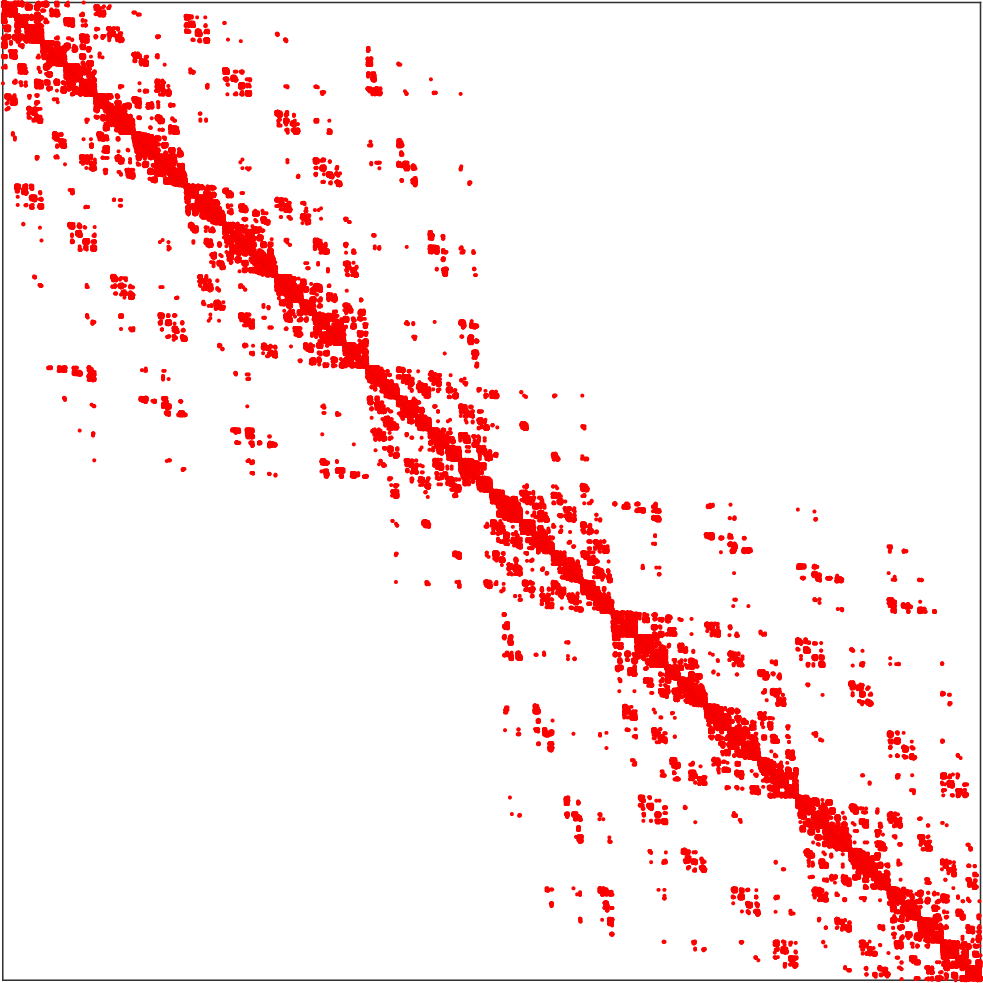}\\
reordered FEM-matrix
\end{minipage}
\begin{minipage}{0.32\textwidth}
\centering
\includegraphics[width=0.95\textwidth]{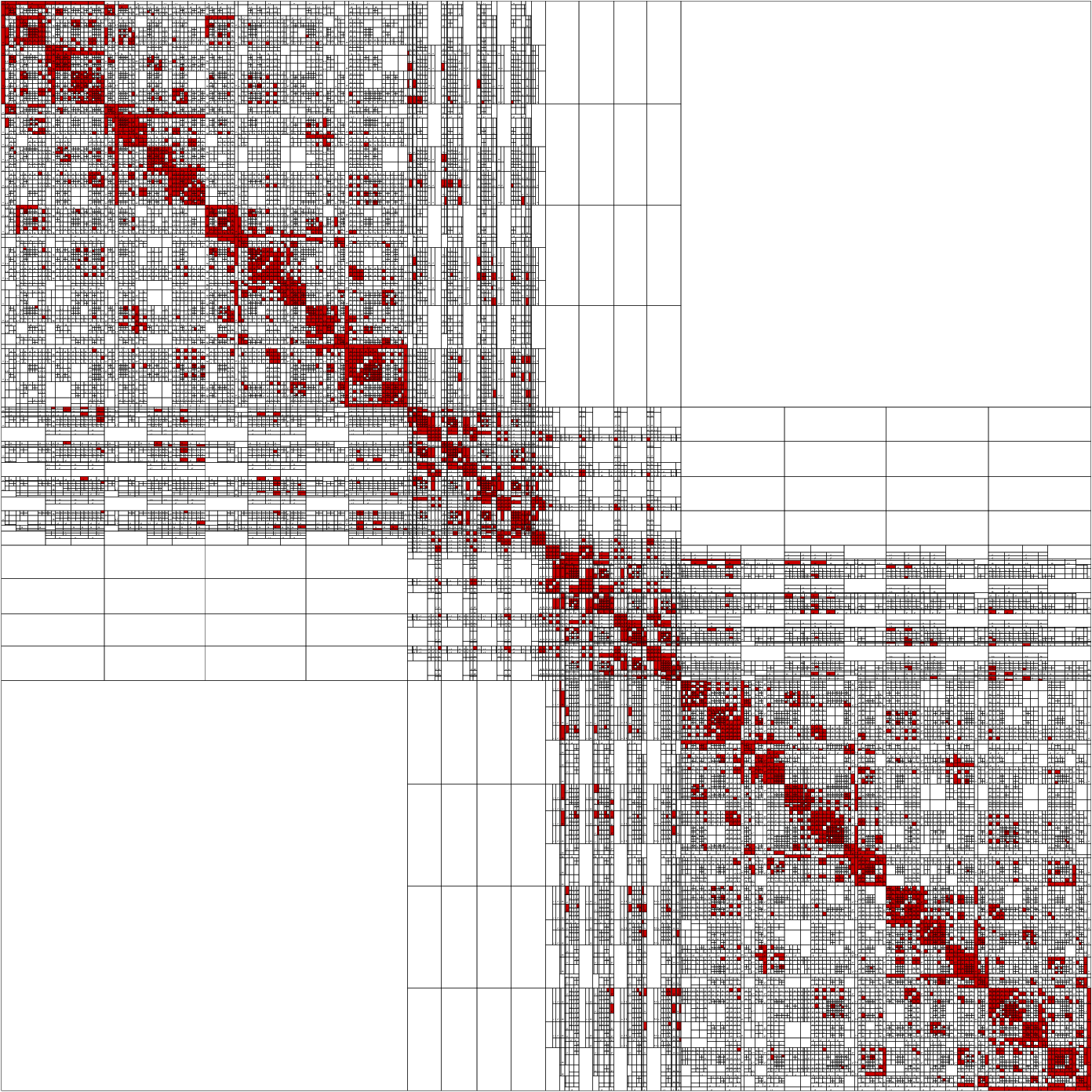}\\
$\mathcal{H}$-matrix representation
\end{minipage}
\caption{\label{fig:FEMHmatrix}Sparsity pattern of a 3D finite element 
matrix,
its reordered finite element matrix, and the corresponding
$\mathcal{H}$-matrix. Red blocks in the $\mathcal{H}$-matrix correspond to the
nearfield, white blocks correspond to the empty farfield.}
\end{figure}

It remains to explain how to actually compute an $\mathcal{H}$-matrix
approximation to the inverse or the LU-decomposition
of a finite element stiffness matrix. To that
end, note that a necessary condition for an entry $\mathbf{A}_{ij}$ in 
the finite element matrix to be non-zero is that $\Upsilon _i\cap\Upsilon _j
\neq\emptyset$, i.e.,\ the intersection of the corresponding supports 
of the basis functions is non-empty. This yields together with the 
$\eta$-adminissibility condition \eqref{eq:admissibility} that all entries 
of a finite element matrix have $\eta$-inadmissible supports, i.e.,\ they 
are contained in the nearfield of an $\mathcal{H}$-matrix. A sparse finite 
element matrix can therefore be represented as an $\mathcal{H}$-matrix 
by reordering the index set corresponding to the clustering scheme introduced 
in Section \ref{sec:ClTr} and inserting the non-zero entries into the nearfield.
An illustration of this procedure can be found in Figure \ref{fig:FEMHmatrix}.

Having the finite element matrix represented by an $\mathcal{H}$-matrix,
the approximate inverse and the LU-decomposition can be computed 
by using the block algorithms of the $\mathcal{H}$-matrix arithmetic 
in $\mathcal{O}(k^2 N\log^2 N)$ operations, cf.~\cite{Beb08,Bor10,Hac15}.
We note especially that the computation of the LU-decomposition 
together with its forward and backward substitution algorithms 
still have an overall complexity of $\mathcal{O}(k^2 N\log^2 N)$, 
but with smaller constants than the computation and the application 
of an approximate inverse.

\subsection{Weak admissibility}
Approximate $\mathcal{H}$-matrix representations for the inverse or
LU-factorizations of finite element matrices have been used to construct
preconditioners for iterative solvers, see, e.g.,~\cite{Beb08} and the
references therein. In \cite{HKK04}, it was observed for the one-dimensional case that the computation
of an approximate inverse can be considerably sped up by replacing the
$\eta$-admissibility condition \eqref{eq:admissibility} by the following
\emph{weak admissibility condition}.

\begin{definition}\label{def:weakadm}
Two clusters $\sigma$ and $\sigma'$ are called \emph{weakly admissible}
if $\sigma\neq\sigma'$.
\end{definition}

We observe immediately that an $\eta$-admissible block cluster is also 
weakly admissible. Thus, by replacing the $\eta$-admissibility condition 
by the weak admissibility, we obtain a much coarser partition of the 
$\mathcal{H}$-matrix. This leads to smaller constants in the storage 
and computational complexity, cf.\ \cite{HKK04}. Each row and each 
column of the finite element matrix has only $\mathcal{O}(1)$ entries. 
Thus, inserting the finite element matrix into a weakly admissible 
$\mathcal{H}$-matrix structure, the off-diagonal blocks of the 
$\mathcal{H}$-matrix have low-rank.

By partitioning the matrix according to the weak admissibility 
condition, we cannot ensure the exponential convergence of fast 
black box low-rank approximation techniques as used for boundary 
element matrices. For example, the adaptive cross approximation (ACA) 
relies on an admissibility condition similar to \eqref{eq:admissibility} to 
ensure exponential convergence, cf.~\cite{Beb00}. Instead, the authors
of \cite{HKK04} suggest to assemble a weakly admissible matrix block 
according to the $\eta$-admissibility condition and transform it on-the-fly 
to a low-rank matrix to obtain a good approximation.

The behavior of the ranks of the low-rank matrices in weakly admissible 
partitions compared to $\eta$-admissible partitions is not fully understood yet.
Suppose that $k_\eta$ is an upper bound for the ranks corresponding to an
$\eta$-admissible partition and suppose that $k_w$ shall be an upper bound for
the ranks to a weakly admissible partition. In \cite{HKK04}, it is proven for
one-dimensional finite element discretizations that one should generally choose
\[
k_w=Lk_\eta,
\]
in order to obtain the same approximation accuracy in the weakly 
admissible case as in the $\eta$-admissible case. Here, $L$ is a 
constant which depends on the depth of the block cluster tree 
and thus logarithmically on $N$. Already in the same article, the 
authors remark in the numerical examples that this bound on 
$k_w$ seems to be too pessimistic and one could possibly choose
\begin{equation}
\label{eq:etaweakconst}
k_w=c_{\eta\rightarrow w}k_\eta,
\end{equation}
where $c_{\eta\rightarrow w}\in[2,3.5]$.

\begin{figure}
\centering
\begin{minipage}{0.4\textwidth}
\centering
\includegraphics[width=0.95\textwidth]{FEM-Hmatrix-eta}\\
$\eta$-admissible FEM-matrix
\end{minipage}
\begin{minipage}{0.4\textwidth}
\centering
\includegraphics[width=0.95\textwidth]{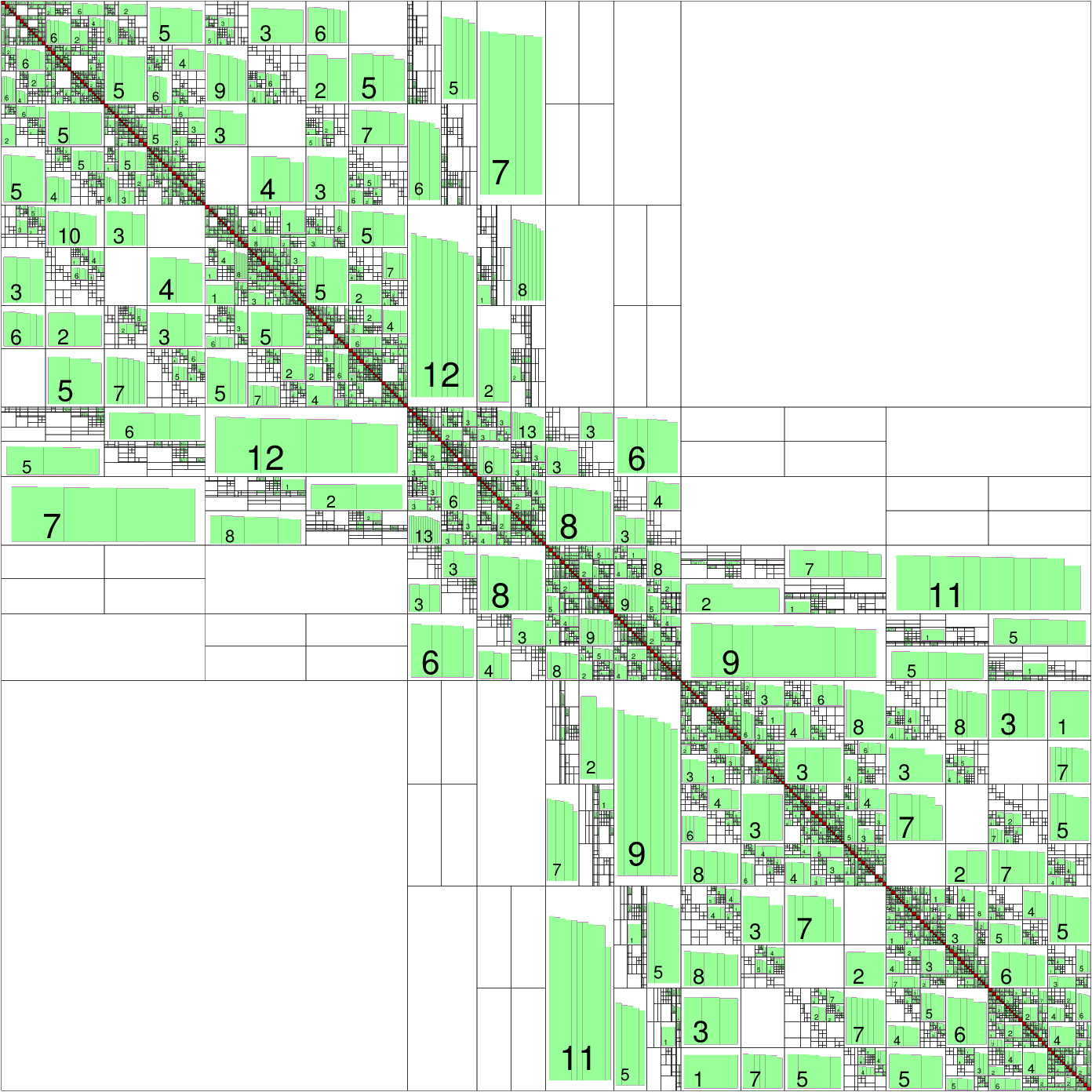}\\
weakly admissible FEM-matrix
\end{minipage}
\caption{\label{fig:weakadm}Comparison of the partition for $\eta$-admissibility
and for weak admissibility. Red blocks correspond to full matrices, green blocks
correspond to low-rank matrices with inscribed rank, and white blocks are zero.}
\end{figure}

Unfortunately, the weak admissibility is not suitable for dimensions greater
than one due to the fact that clusters can possibly intersect each other in
$\mathcal{O}(N^\alpha)$ points, where $\alpha\ge 0$ depends on the spatial
dimension. However, one can try to reduce this negative influence of the weak
admissibility condition by mixing it with the $\eta$-admissibility.
In the software package HLib, cf.~\cite{BGH03}, 
the authors use the $\eta$-admissibility for all block clusters with a block size 
larger than a given threshold and apply the weak admissibility condition for block 
clusters $\sigma\times\sigma'$ which are below that threshold provided that
the condition
\[a_i^\sigma<\frac{a_i^{\sigma '}+b_i^{\sigma '}}{2}<b_i^\sigma\quad\text{or}\quad
a_i^{\sigma '}<\frac{a_i^\sigma+b_i^\sigma}{2}<b_i^{\sigma '}\]
is satisfied for
the corresponding bounding boxes $\prod _{i=1}^3[a_i^\mu,b_i^\mu]$,
$\mu=\sigma,{\sigma '}$, in at most one coordinate direction. This
condition restricts the application of the weak admissibility to essentially one-dimensional cluster intersections with length below a certain threshold. The
impact of this specific admissibility condition is illustrated in
Figure~\ref{fig:weakadm}.

\subsection{Nested dissection}
\label{sec:ND}
While the weak admissibility takes the sparsity of the finite element
matrix into account only during the construction of the
$\mathcal{H}$-matrix, it is also possible to incorporate the sparsity
already during the construction of the cluster tree. One possibility
to do so was introduced in \cite{GKL09} and is based on nested
dissection, cf.~\cite{BT02,Geo73,HR98,LRT79} and the references therein.
We briefly review the idea of nested dissection in the context of 
$\mathcal{H}$-matrices as discussed in \cite{GKL09} and refer 
to \cite{GKL09} for more details.

\begin{figure}
\centering
\begin{minipage}{0.32\textwidth}
\centering
\includegraphics[width=0.95\textwidth]{FEM}\\
sparse FEM-matrix\\
\mbox{}
\end{minipage}
\begin{minipage}{0.32\textwidth}
\centering
\includegraphics[width=0.95\textwidth]{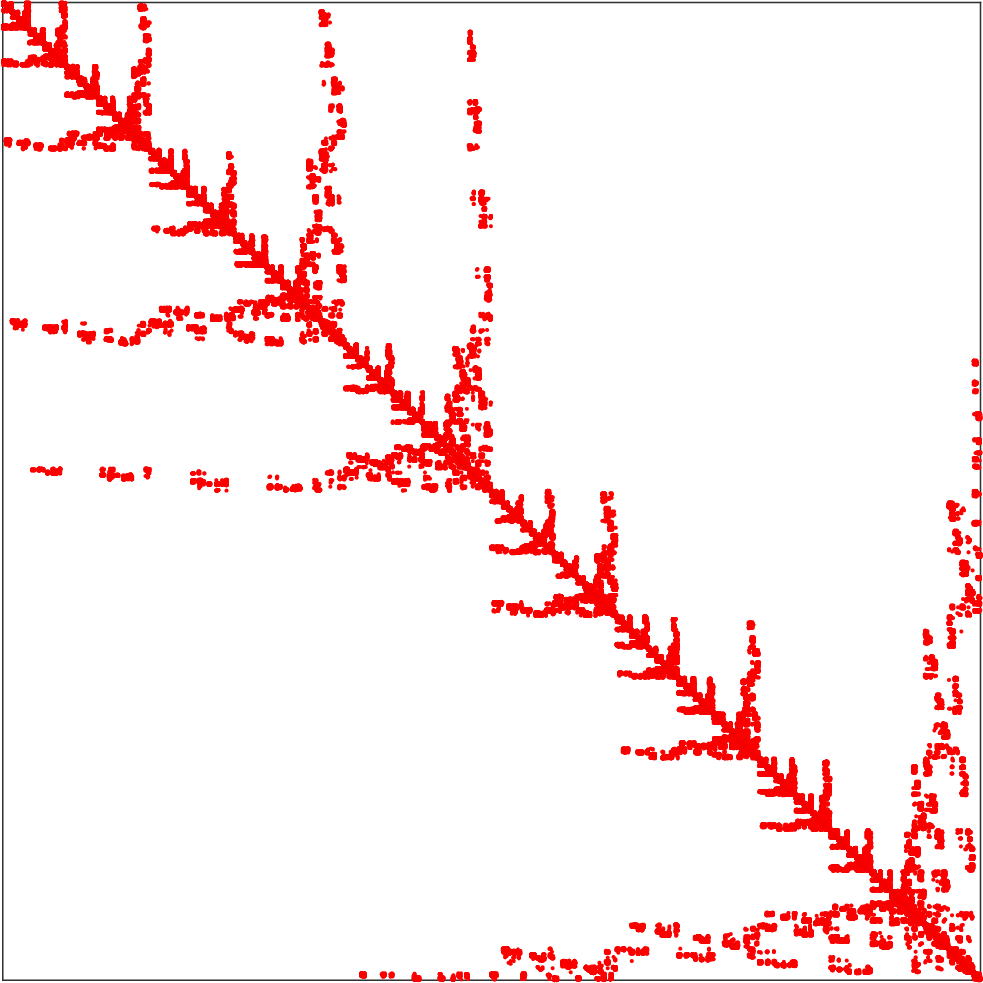}\\
nested dissection\\
reordered FEM-matrix
\end{minipage}
\begin{minipage}{0.32\textwidth}
\centering
\includegraphics[width=0.95\textwidth]{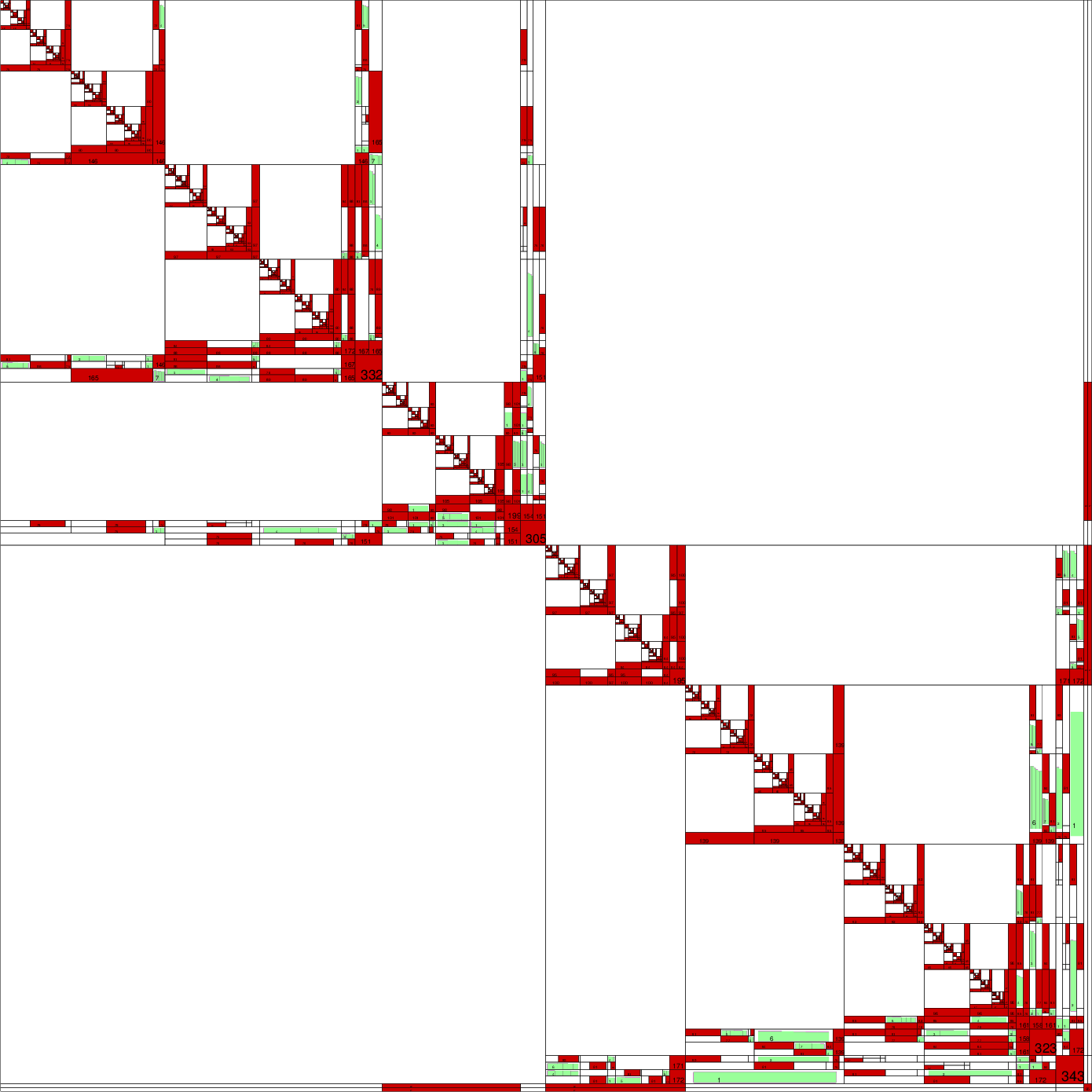}\\
nested dissection\\
$\mathcal{H}$-matrix representation\\
\end{minipage}
\caption{\label{fig:DDHmatrix}Sparsity pattern of a 3D finite element 
matrix, its reordered finite element matrix according to nested
dissection and its corresponding $\mathcal{H}$-matrix.}
\end{figure}

The idea is to employ a recursive algorithm as follows.
\begin{enumerate}
\item Split degrees of freedom into three disjoint subsets
$\mathcal{I}_1$, $\mathcal{I}_2$, $\mathcal{I}_3$ according 
to the following conditions.
\begin{itemize}
\item $\mathcal{I}_1$ and $\mathcal{I}_2$ should have comparable sizes.
\item $\mathcal{I}_1$ and $\mathcal{I}_2$ should not interact with
each other, i.e., all entries $\mathbf{A}_{i,j}$, $i\in\mathcal{I}_1$,
$j\in\mathcal{I}_2$ of the finite element matrix are zero.
\item $\mathcal{I}_3$ is the boundary layer between $\mathcal{I}_1$
and $\mathcal{I}_2$.
\end{itemize}
\item Relabel indices in subsequent order: first $\mathcal{I}_1$, then
$\mathcal{I}_2$, and then $\mathcal{I}_3$.
\item Proceed recursively with $\mathcal{I}_1$ and $\mathcal{I}_2$
(in the $\mathcal{H}$-matrix framework, $\mathcal{I}_3$ will also receive a recursive treatment).
\end{enumerate}
Reordering the index sets of the finite element matrix in accordance
with this procedure yields a sparsity pattern as illustrated in
Figure~\ref{fig:DDHmatrix}. Due to the special construction, 
large parts of the matrix are zero and will remain zero in a 
subsequent LU-decomposition.

In the following, we take the approach of \cite{GKL09} to construct
an $\mathcal{H}$-matrix which reorders the index set such that 
the pattern of the finite element matrix exposes a nested dissection 
ordering. We therefore recapitulate the construction of a cluster tree 
based on domain decomposition as proposed in \cite{GKL09}. For 
that purpose, the cluster algorithm distinguishes between \emph{domain
clusters} and \emph{interface clusters}. The following algorithm is 
employed for the root $\{1,\ldots,N\}$ and all domain clusters.
\begin{enumerate}
\item Given a cluster $\sigma$ and its corresponding set of points
$\{\mathbf{x}_1,\ldots,\mathbf{x}_{\#\sigma}\}$, construct an 
axis-parallel bounding box $Q_\sigma$.
\item Cut $Q_\sigma$ into two pieces $Q_1$ and $Q_2$ by
halving the longest edge.
\item Define three disjoint sons of $\sigma$ as follows.
\begin{itemize}
\item $\sigma_1 = \{i\colon\mathbf{x}_i\in Q_1\}$,
\item $\sigma_2 = \{i\colon\Upsilon_i\cap \Upsilon_{\sigma _1}=\emptyset\}$,
\item $\sigma_3 = \sigma\setminus\{\sigma_1\cup\sigma_2\}$.
\end{itemize}
\item Relabel indices of clusters in subsequent order: first $\sigma_1$, then
$\sigma_2$, and then $\sigma_3$.
\item Treat $\sigma_1$ and $\sigma_2$ as domain clusters and $\sigma_3$ as
interface cluster as indicated below.
\end{enumerate}
Due to their special construction, the bounding boxes of interface clusters are
``flat'' in one coordinate direction. The cluster algorithm has therefore to be
adapted to fit the asymptotic requirements of Section~\ref{sec:ClTr}. Therefore,
we define $\operatorname{level}_{\text{int}}(\sigma)$ as the distance of
$\sigma$ to the nearest domain cluster in the cluster tree. We then employ
the following cluster algorithm for the interface clusters and refer to
\cite{GKL09} for a detailed discussion.
\begin{enumerate}
\item Distinguish two cases:
\begin{itemize}
\item If $\operatorname{level}_{\text{int}}(\sigma)=0\mod d$: do not subdivide
$\sigma$ and set $\sigma'=\sigma$ as its only son.
\item If $\operatorname{level}_{\text{int}}(\sigma)\neq 0\mod d$: split the
bounding box $Q_\sigma$ axis parallel in two boxes $Q_1$ and $Q_2$ such that
the ``flat'' direction is not modified. Then define the two son clusters
$\sigma _1=\{i\colon \mathbf{x}_i\in Q_1\}$ and
$\sigma _2=\{i\colon \mathbf{x}_i\in Q_2\}$.
\end{itemize}
\item Apply the cluster algorithm for interface clusters to all sons.
\end{enumerate}

In order to translate the sparsity of the finite element matrix into the
block structure of an $\mathcal{H}$-matrix, we can combine the
$\eta$-admissibility condition from \eqref{eq:admissibility} and the
weak admissibility condition from Definition~\ref{def:weakadm} to a
\emph{nested dissection admissibility} condition.
\begin{definition} \label{def:ndadm}
Two clusters $\sigma$ and $\sigma'$ are called \emph{nd-admissible}
if either
\begin{itemize}
\item $\sigma\neq\sigma '$ are both domain clusters or
\item $\sigma$ and $\sigma'$ are $\eta$-admissible.
\end{itemize}
\end{definition}
In fact, if $\sigma\neq\sigma '$ are both domain clusters, we can directly say
that the corresponding $\mathcal{H}$-matrix block has rank zero.
Figure~\ref{fig:DDHmatrix} illustrates the sparsity pattern of the finite
element matrix after the permutations determined by the cluster algorithms
and how large parts of the corresponding $\mathcal{H}$-matrix have rank zero.
The low-rank blocks in the representation are due to some internal checks of
HLib, which aim at replacing inadmissible blocks by low-rank matrices only if
very few entries in the corresponding matrix block are non-zero.

The numerical experiments in the next section show that the sparse 
structure constructed here leads to smaller constants in the complexity
of the solution algorithm.

\section{Numerical results}\label{sec:numres}
Before we summarize the settings of the numerical experiments, we briefly
recall that the algorithm of the presented method consists of the following
three steps.
\begin{enumerate}
\item Compute the sparse finite element matrix $\mathbf{A}$ in linear and the
correlation $\mathcal{H}$-matrix $\mathbf{C}_f$ in almost linear time.
\item Compute the approximate LU-decomposition of $\mathbf{A}$ in
$\mathcal{H}$-matrix format in almost linear time.
\item Solve the matrix equation \eqref{eq:linsys2} with iterative refinement
\eqref{eq:iterativerefinement} in almost linear time for each iteration.
\end{enumerate}
The numerical experiments in this article shall mainly focus on the third step
and the overall behaviour of the method. We will see that only one iteration is
required in the third step, which yields an almost linear overall complexity of
the algorithm.
To improve the computation time, we store and compute only the 
lower part of $\mathbf{C}_f$ and $\mathbf{C}_u$.

All the computations in the following experiments have been carried 
out on a single core of a computing server with two Intel(R) Xeon(R) 
E5-2670 CPUs with a clock rate of 2.60GHz and a main memory of 
256GB. For the $\mathcal{H}$-matrix computations, we use the software
package HLib, cf.~\cite{BGH03}, and for the finite element discretization 
the Partial Differential Equation Toolbox of {\sc Matlab}\footnote{Release 2015b, 
The MathWorks, Inc., Natick, Massachusetts, United States.} which employs
piecewise linear finite elements. The two libraries are coupled together 
in a single C-program, cf.~\cite{KR88}, using the {\sc Matlab} Engine 
interface. The meshes are generated by {\sc Tetgen}, cf.~\cite{Si15},
and then imported into {\sc Matlab}.

\subsection{Experimental setup}
To obtain computational efficiency and to keep the ranks of the low-rank 
matrices under control, HLib imposes an upper threshold $k_\eta$ for the 
ranks in the case of an $\eta$-admissible $\mathcal{H}$-matrix and a lower 
threshold $n_{\min}$ for the minimal block size. For the application of the weak 
admissibility condition, we rely on the criterion of HLib, which considers the 
weak admissibility condition only if one of the index sets of the block cluster 
$\sigma\times\sigma'$ has a cardinality below 1,024 and
the condition 
\[a_i^\sigma<\frac{a_i^{\sigma '}+b_i^{\sigma '}}{2}<b_i^\sigma\quad\text{or}\quad
a_i^{\sigma '}<\frac{a_i^\sigma+b_i^\sigma}{2}<b_i^{\sigma '}\] 
is satisfied for
the corresponding bounding boxes $\prod _{i=1}^3[a_i^\mu,b_i^\mu]$,
$\mu=\sigma,{\sigma '}$, in at most one coordinate direction. Otherwise,
the $\eta$-admissibility is used instead. In the case of a weakly admissible 
matrix block, HLib imposes an upper threshold of $k_w=3k_\eta$, setting 
$c_{\eta\rightarrow w}=3$ in \eqref{eq:etaweakconst}. For our experiments, 
we choose $\eta=2$, $k_{\eta}=20$, $n_{\min}=50$ and employ either a 
geometric cluster strategy, i.e.,~the binary cluster strategy from the end of 
Section~\ref{sec:ClTr} or the nested dissection cluster strategy
as discussed in Section~\ref{sec:ND}.
The iterative refinement is stopped if the absolute 
error of the residual in the Frobenius norm is smaller than $10^{-6}$.

In the following examples, we want, besides other aspects, to study
the influence of the weak admissibility condition and the nested
dissection clustering for the partitioning of the
different $\mathcal{H}$-matrices. Namely, we successively want to replace 
the $\eta$-admissibility by the weak admissibility for a binary and a
nested dissection cluster tree as described in Table
\ref{tab:admissibility} in order to lower the constants hidden in the complexity 
of the $\mathcal{H}$-matrix arithmetic and thus to improve the computation 
time. For the discretization of the correlation kernel $\Cor _f$, we will 
always use ACA.

\begin{table}[h]
\centering
\begin{center}
\begin{tabular}{|c|c|c|c|c|}
\hline
\multirow{2}{*}{Case} & \multicolumn{2}{|c|}{Operator and admissibility} & 
\multirow{2}{*}{Cluster tree} \\\cline{2-3}
           & $\mathcal{L}$ and $\mathcal{L}^{-1}$
                                  & $\Cor _f$ and $\Cor _u$ & \\\hline\hline
all-$\eta$ & $\eta$-admissibility & $\eta$-admissibility    & binary \\\hline
weak-FEM   & weak admissibility   & $\eta$-admissibility    & binary \\\hline
all-weak   & weak admissibility   & weak admissibility      & binary \\\hline
nd-$\eta$  & nd-admissibility     & $\eta$-admissibility & nested dissection    \\\hline
nd-weak    & nd-admissibility     & weak-admissibility   & nested dissection    \\\hline
\end{tabular}
\end{center}
\caption{\label{tab:admissibility}The five combinations of the 
	admissibility conditions used for the partition of the
         $\mathcal{H}$-matrices.}
\end{table}

The \emph{all-$\eta$} case is the canonical case and has also been 
investigated in case of the boundary element method in \cite{DHP15}. 
The \emph{weak-FEM} case is a first relaxation to apply the weak admissibility
condition. This is justified, since the stiffness matrix $\mathbf{A}$ can exactly
be represented as a weakly admissible $\mathcal{H}$-matrix and the iterative
refinement only involves an approximate LU-decomposition.
Hence, we expect at most an influence on the quality
of the approximate LU-decomposition and thus on the number of iterations
in the iterative refinement. We have therefore to investigate if possible
additional iterations are compensated by the faster $\mathcal{H}$-matrix
arithmetic.

The aforementioned cases have in common that they rely on the 
asymptotic smoothness of $\Cor _f$ and $\Cor_u$ and the $\eta$-admissibility 
which leads to exponential convergence of the $\mathcal{H}$-matrix approximation.
In the case \emph{all-weak}, we want to examine if there is some indication that 
the weak admissibility could possibly also be considered for the partition of the 
$\mathcal{H}$-matrices for $\Cor _f$ and $\Cor _u$. To that end, we approximate 
$\Cor _f$ with ACA relative to the $\eta$-admissibility partition and convert it 
on-the-fly to the partition of the weak admissibility, as proposed in \cite{HKK04}.

While the three aforementioned cases all rely on a binary cluster tree, the
cases \emph{nd-$\eta$} and \emph{nd-weak} rely on a cluster tree which is 
constructed by nested dissection. In both cases, the finite element matrix 
$\mathbf{A}$ is partitioned by the nd-admissibility. For $\Cor_f$ and $\Cor_u$,
we use the $\eta$-admissibility in the nd-$\eta$-case and the weak admissibility 
in the nd-weak-case, where we assemble the matrix for $\Cor_f$ in the same way 
as in the all-weak-case.

The following numerical examples are divided in two parts. In the first part,
we demonstrate the convergence of the presented method by comparing it to
a low-rank reference solution computed with the pivoted Cholesky decomposition,
cf.~\cite{HPS12}. In the second part, we will demonstrate that the presented
method works also well in the case of correlation kernels with low Sobolev
smoothness or small correlation length, where no low-rank approximations
exist and sparse tensor product approximations fail to resolve the correlation length.
Note that in both examples, due to the non-locality of the
correlation kernels and the Green's function, the computed system matrices
are smaller than usual for the finite element method. In particular, the unknown
in the system of equations \eqref{eq:linsys2} is a matrix with
$N^2$ entries, whereas the corresponding mesh has only $N$ degrees of freedom.
The $\mathcal{H}$-matrix compression reduces the computational complexity for
the assembly and the amount of required
storage from $N^2$ to $\mathcal{O}(kN\log N)$, whereas the complexity of the
solution algorithm decreases from $\mathcal{O}(N^3)$ to $\mathcal{O}(k^2N\log^2 N)$.

\subsection{Tests for the iterative solver}
Due to the recompression schemes in the block matrix algorithms 
of the $\mathcal{H}$-matrix arithmetic, it is not directly clear if the 
presented solver converges. Still, it can be shown that 
some iterative $\mathcal{H}$-matrix schemes converge up to a certain 
accuracy, cf.~\cite{HKT08}. In the following, we want to demonstrate 
for a specific example that our iterative scheme provides indeed 
convergence.

\begin{figure}[hbt]
\centering
\includegraphics[height=0.23\textwidth, clip=true, trim=120 530 140 110]{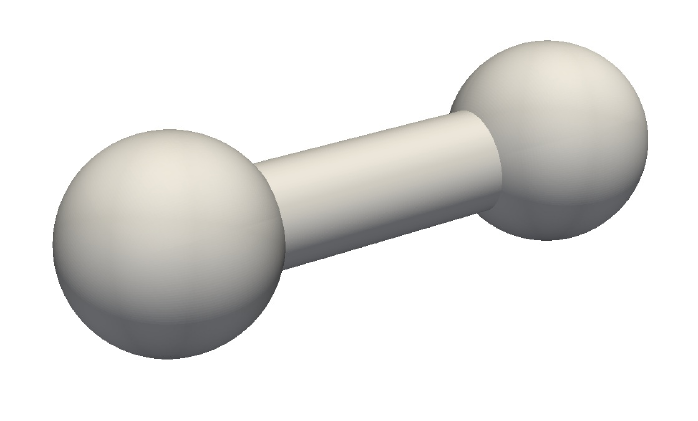}\quad
\includegraphics[height=0.23\textwidth, clip=true, trim=120 530 30 110]{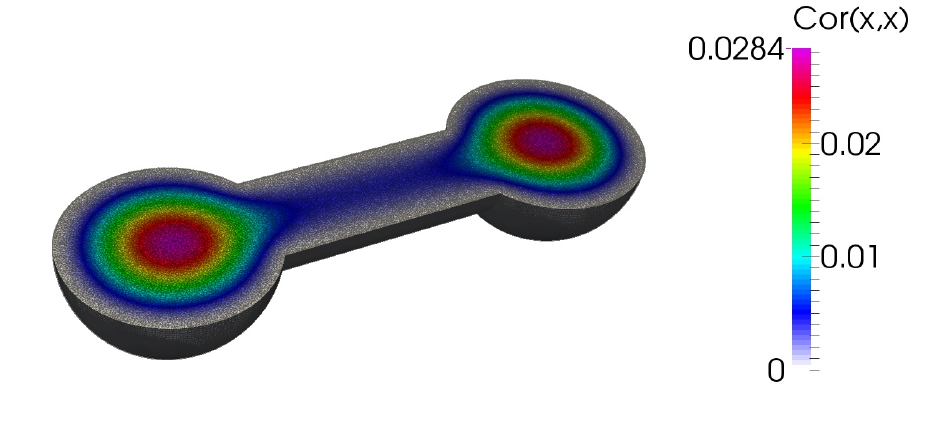}
\caption{\label{fig:bone}The dumbbell geometry and its meshed cross section
with the trace of the solution correlation kernel $\Cor _u|
_{\mathbf{x}=\mathbf{y}}$ for load data prescribed by the Mat\'ern-5/2 kernel.}
\end{figure}

On the dumbbell geometry pictured in Figure~\ref{fig:bone}, we 
consider $\mathcal{L} = -\Delta$ in \eqref{eq:coreq} and the
Mat\'ern-5/2 kernel as input correlation $\Cor _f$, i.e.,\ for 
$r=\|\mathbf{x}-\mathbf{y}\|_2$, we set
\[
\Cor_f(\mathbf{x},\mathbf{y})=\bigg(1+\sqrt{5}\frac{r}{\ell}
+\frac{5}{3}\frac{r^2}{\ell^2}\bigg)\exp\Big(-\sqrt{5}\frac{r}{\ell}\Big),
\]
where $\ell\approx\diam (D)$ denotes the correlation length. The 
conversion of the finite element matrix to an $\mathcal{H}$-matrix 
for the dumbbell geometry has already been illustrated in Figure
\ref{fig:FEMHmatrix}. Whereas, the difference between the 
$\eta$-admissibility and the weak admissibility is illustrated
in Figure~\ref{fig:weakadm} and the effect of the nested dissection
ordering is illustrated in Figure~\ref{fig:DDHmatrix}.

\begin{table}[hbt]
\centering
\scalebox{0.98}{
\begin{tabular}{|c|cccccc|c|}
\hline
Level & 1 & 2 & 3 & 4 & 5 & Reference mesh \\\hline\hline
Mesh points  &        238 &      1,498 &      6,958 &     34,112 &    175,562 &   1,033,382\\\hline
$N$          &          4 &        201 &      1,742 &     13,341 &     98,177 &   756,626\\\hline
$N^2$        &      $16$ &    40,401 &  3,034,564 & $1.78\cdot10^8$ 
& $9.64\cdot10^9$ &$5.72\cdot10^{11}$\\\hline
\end{tabular}}
\caption{\label{tab:bonedof}Mesh points and degrees of freedom of the
         finite element mesh $N$ and number of entries $N^2$ in
         the solution matrix for different levels of the dumbbell
         geometry.}
\end{table}

For determining a reference solution, we compute a low-rank 
approximation $\mathbf{C}_f\approx\mathbf{L}_f\mathbf{L}_f^{\intercal}$
with the pivoted Cholesky decomposition as proposed in \cite{HPS12}.
The numerical solution $\mathbf{C}_u$ of \eqref{eq:linsys2} is then given by
\[
\mathbf{C}_u \approx \mathbf{L}_u\mathbf{L}_u^{\intercal},
\]
where $\mathbf{L}_u$ solves $\mathbf{A}\mathbf{L}_u=\mathbf{L}_f$.
To compute the error of the $\mathcal{H}$-matrix solution, we compare 
the correlation and the correlation's trace $\Cor _{u,N}|_{\mathbf{x}=\mathbf{y}}$
of the  $\mathcal{H}$-matrix approximation with the correlation
and the correlation's trace $\Cor _{u,N}|_{\mathbf{x}
=\mathbf{y}}$ derived from the pivoted Cholesky decomposition on a 
finer reference mesh. We refer to Table \ref{tab:bonedof} for more 
details on the meshes under consideration.

\begin{figure}[hbt]
\centering
\includegraphics[width=0.35\textwidth]{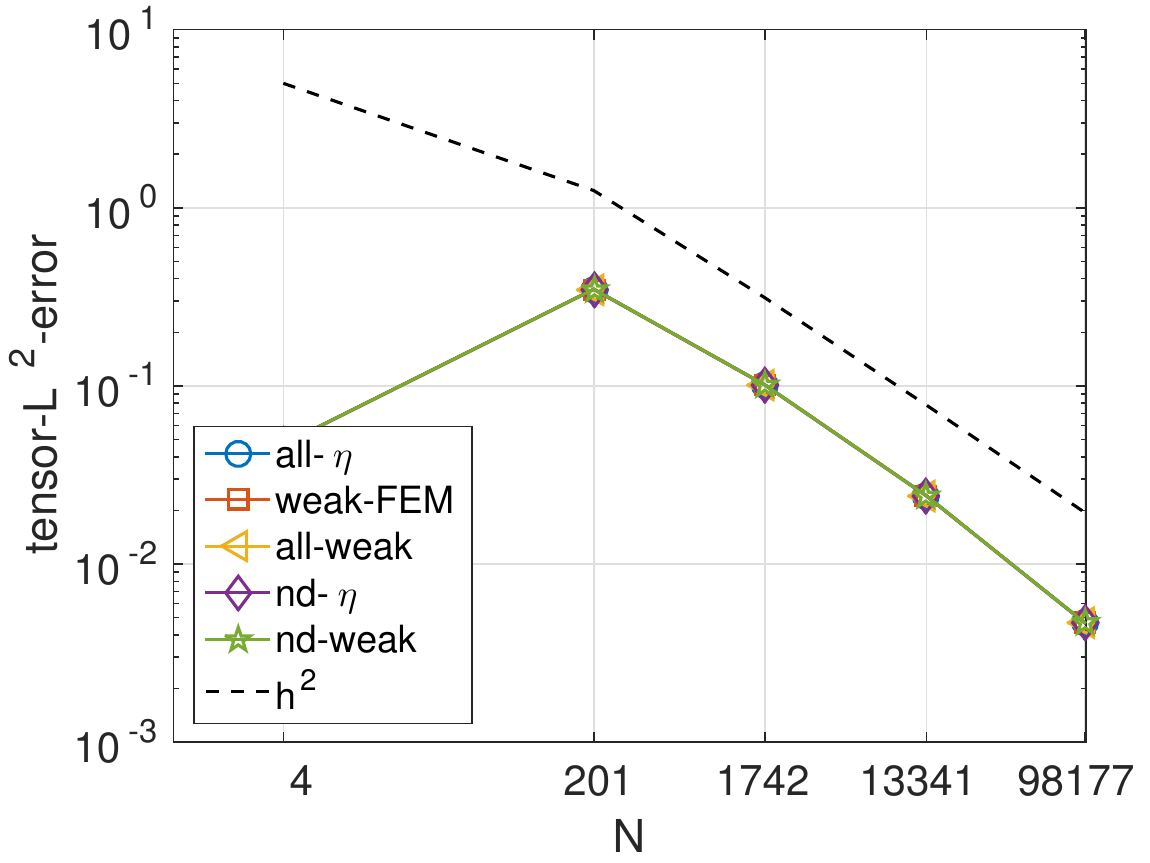}\quad
\includegraphics[width=0.35\textwidth]{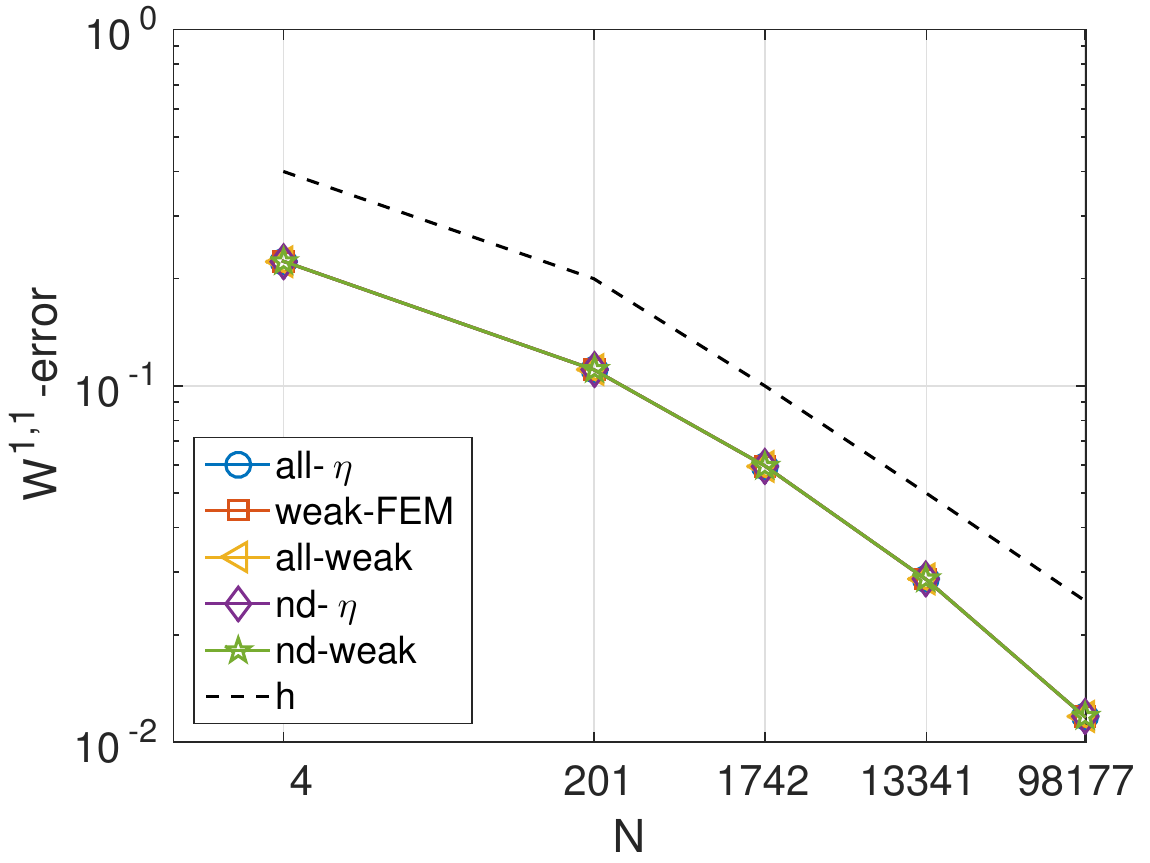}\quad
\caption{\label{fig:boneconvergence}$L^2$-error of $\Cor_{u,N}$ (left)
         and $W^{1,1}$-error of $\Cor _{u,N}|_{\mathbf{x}=\mathbf{y}}$ (right).}
\end{figure}

While the error of the correlation itself can be measured in the $L^2$-norm
on the tensor product domain, the appropriate norm for error measurements of
it's trace is the $W^{1,1}$-norm.
Due to the Poincar\'e-Friedrich and the Cauchy-Schwartz inequality
and
\[
\big\| u^2-u_h^2\big\|_{W^{1,1}(D)}
\lesssim \big\|\nabla \big(u^2-u_h^2\big)\big\|_{L^1(D)}
\lesssim |u-u_h|_{H^1(D)}+\|u-u_h\|_{L^2(D)}
\lesssim h,
\]
we can expect a convergence rate in the $W^{1,1}$-norm 
which is proportional to the mesh size $h$. A standard tensor 
product argument yields a convergence rate of order $h^2$ in
the tensor product $L^2$-norm. Figure~\ref{fig:boneconvergence} 
shows that we indeed reach these rate for all five cases of admissibility 
which are considered in Table \ref{tab:admissibility}. In fact, the observed 
errors coincide in the first few digits.

\begin{figure}[hbt]
\centering
\includegraphics[width=0.35\textwidth]{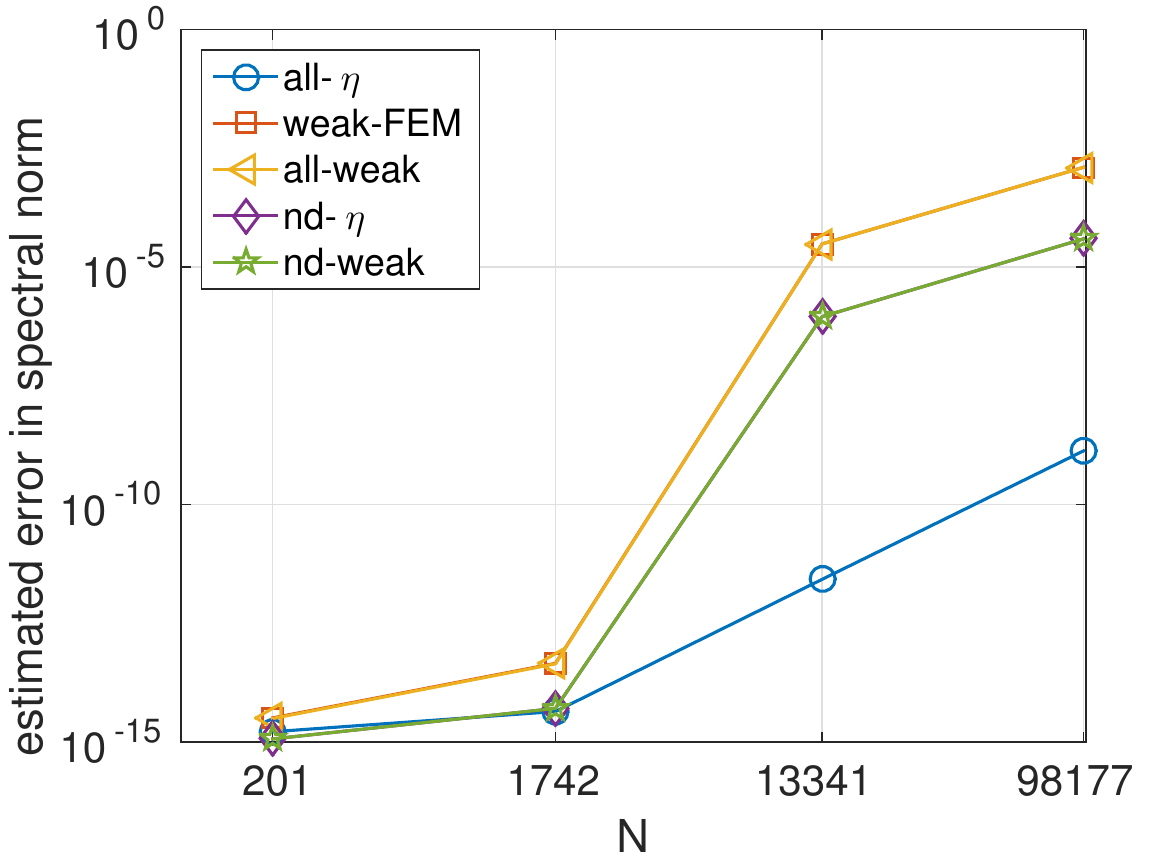}
\caption{\label{fig:deviation}The deviation $\hat{\mathbf{L}}\hat{\mathbf{U}}
         -\mathbf{A}$ in the estimated spectral norm for fixed rank.}
\end{figure}

\begin{figure}[hbt]
\centering
\includegraphics[width=0.3\textwidth]{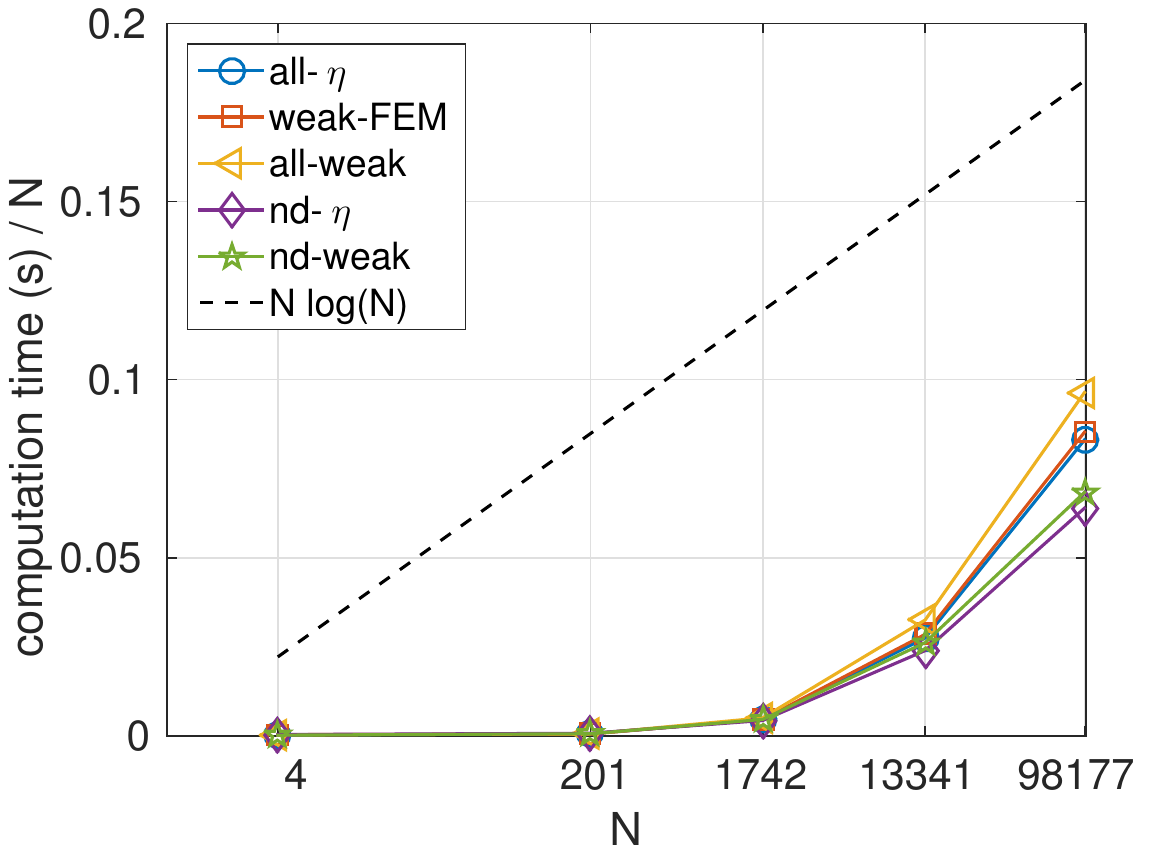}\quad
\includegraphics[width=0.3\textwidth]{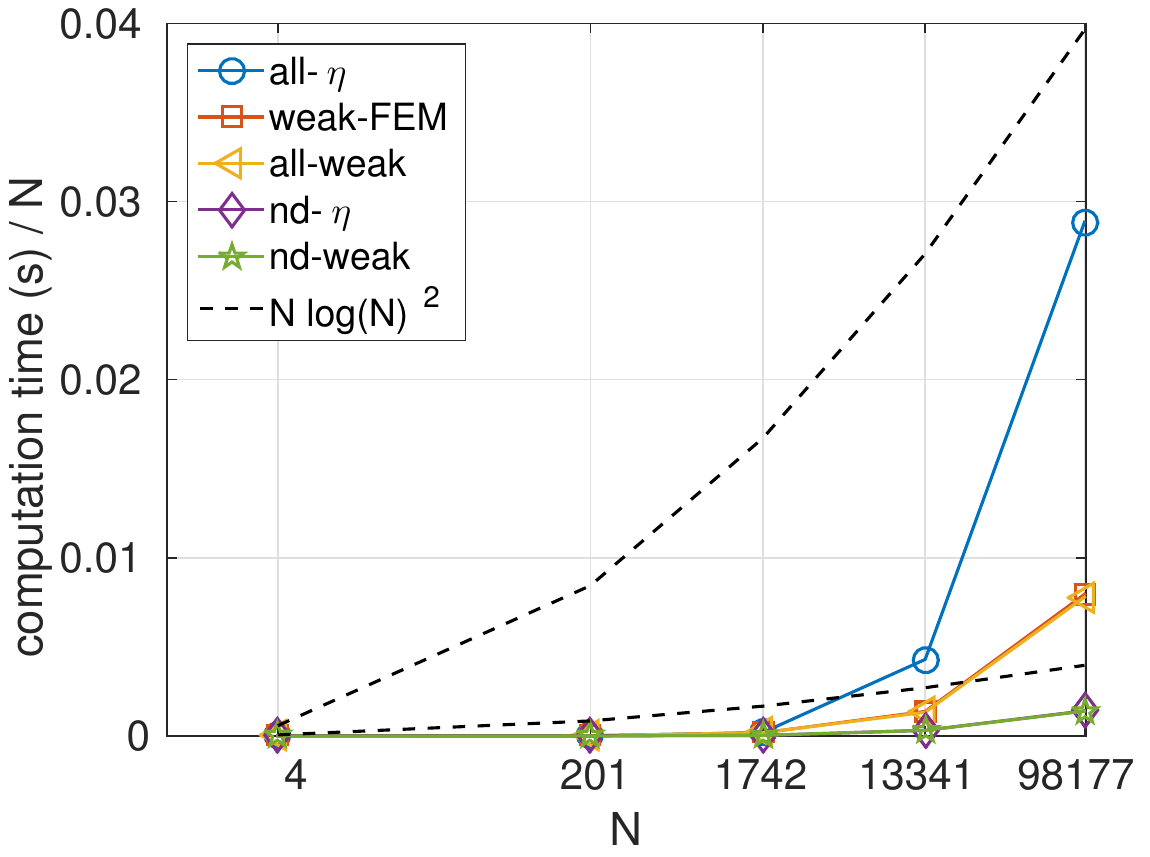}\quad
\includegraphics[width=0.3\textwidth]{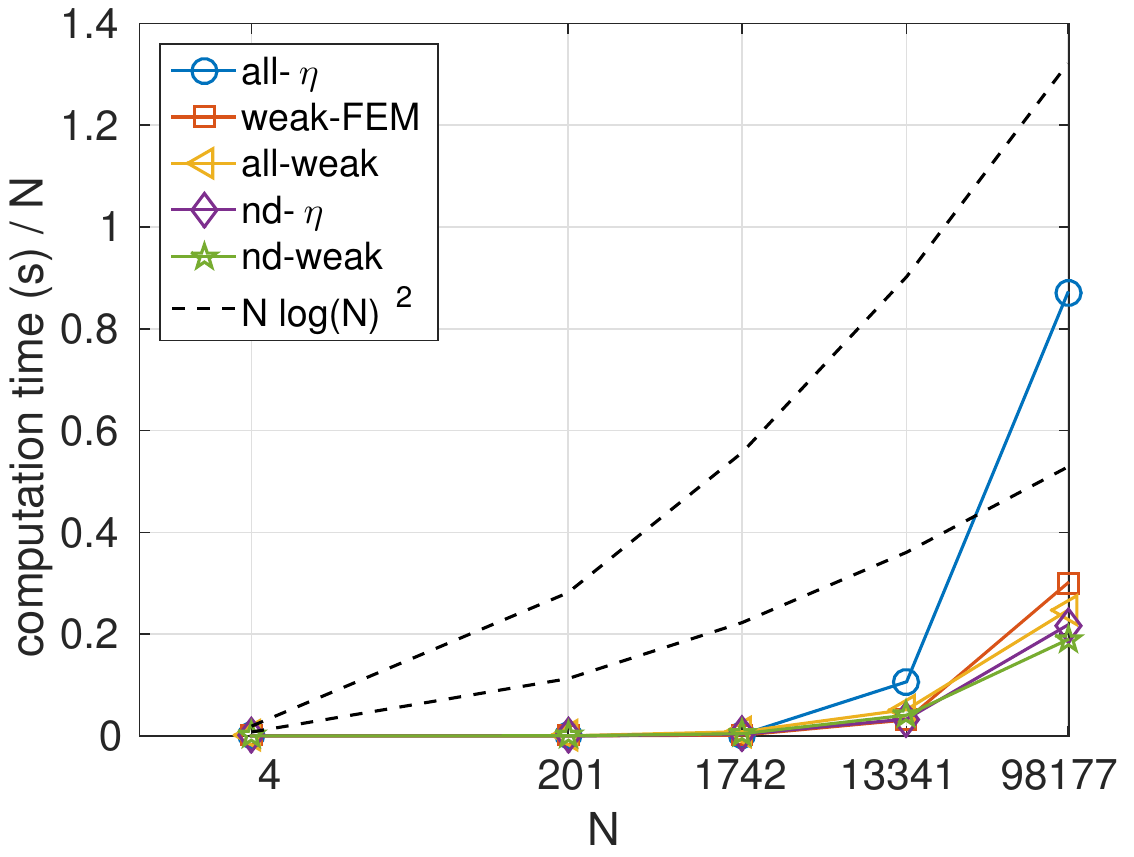}
\caption{\label{fig:bonetimes}Computation times in seconds for the computation
of the data correlation $\mathcal{H}$-matrix $\mathbf{C}_f$ (left), the
approximate LU-decomposition of the system matrix $\mathbf{A}$ (center) and for
the iterative refinement (right) on the dumbbell geometry.}
\end{figure}

\begin{figure}[hbt]
\centering
\includegraphics[width=0.3\textwidth]{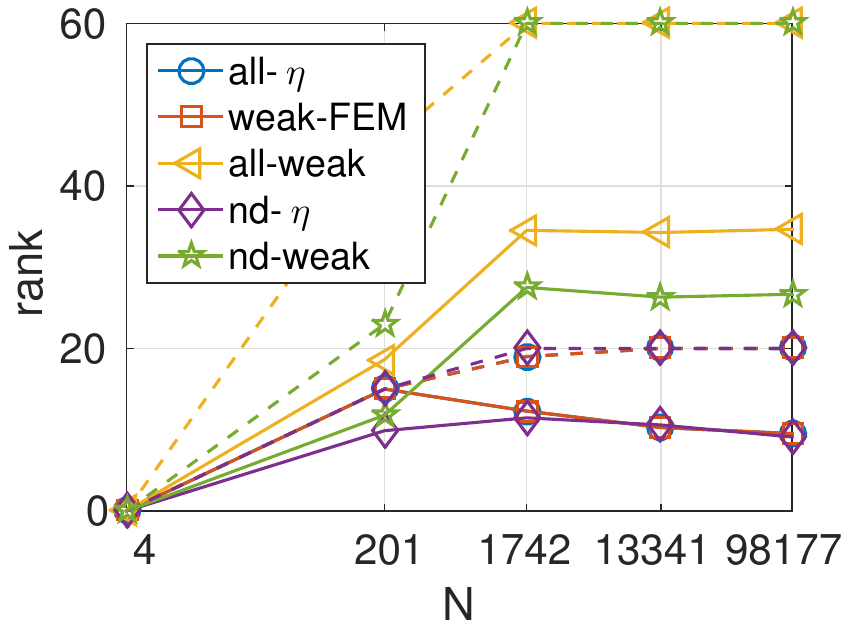}\quad
\includegraphics[width=0.3\textwidth]{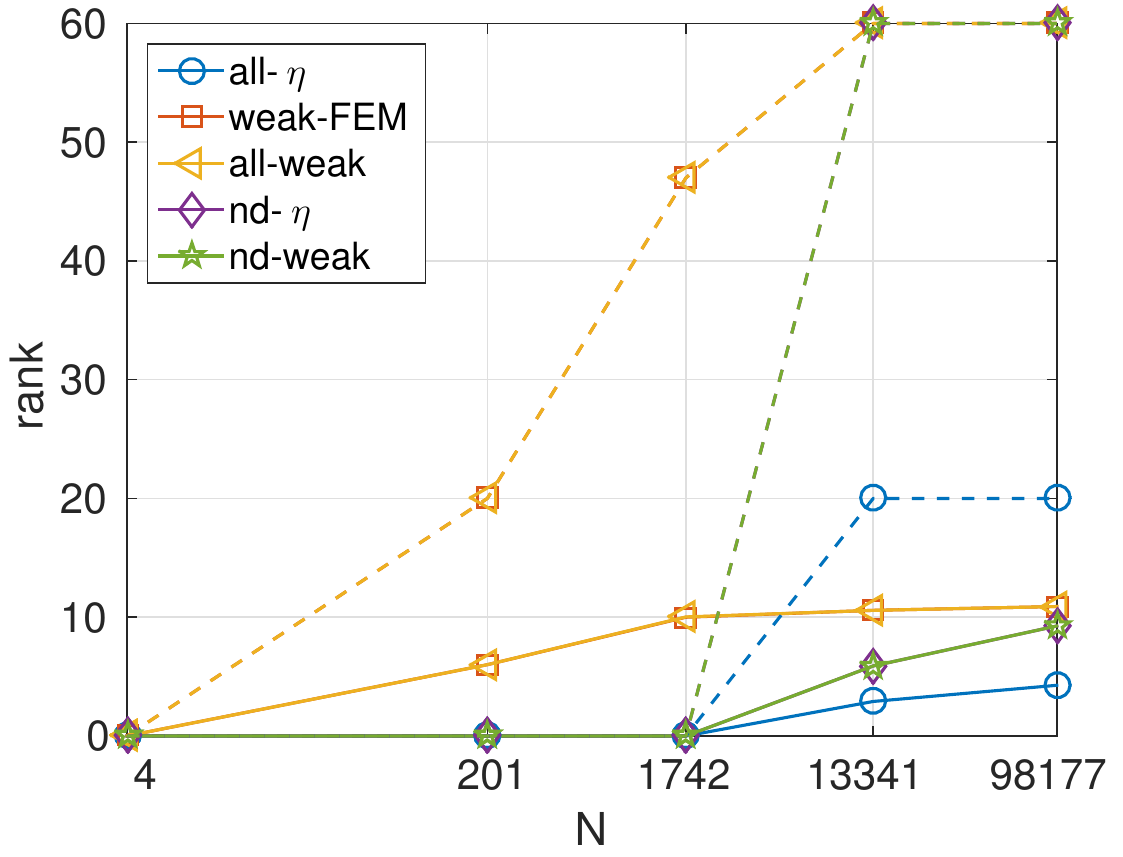}\quad
\includegraphics[width=0.3\textwidth]{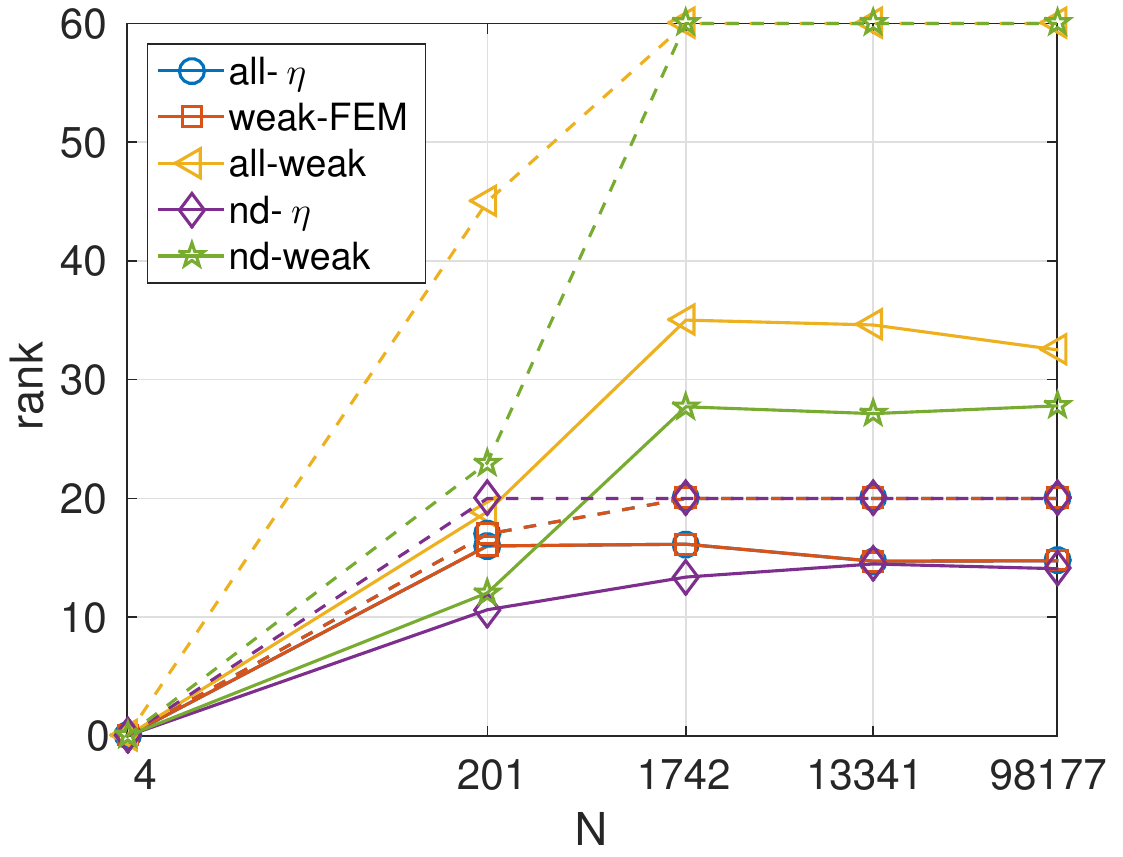}
\caption{\label{fig:boneranks}Required ranks for the prescribed
correlation $\mathbf{C}_f$ (left), the LU-decomposition of $\mathbf{A}$ (middle)
and the solution correlation $\mathbf{C}_u$
(right). The straight lines indicate the average ranks, whereas the dashed
lines illustrate the maximal rank attained.}
\end{figure}

\begin{figure}[hbt]
\centering
\includegraphics[width=0.35\textwidth]{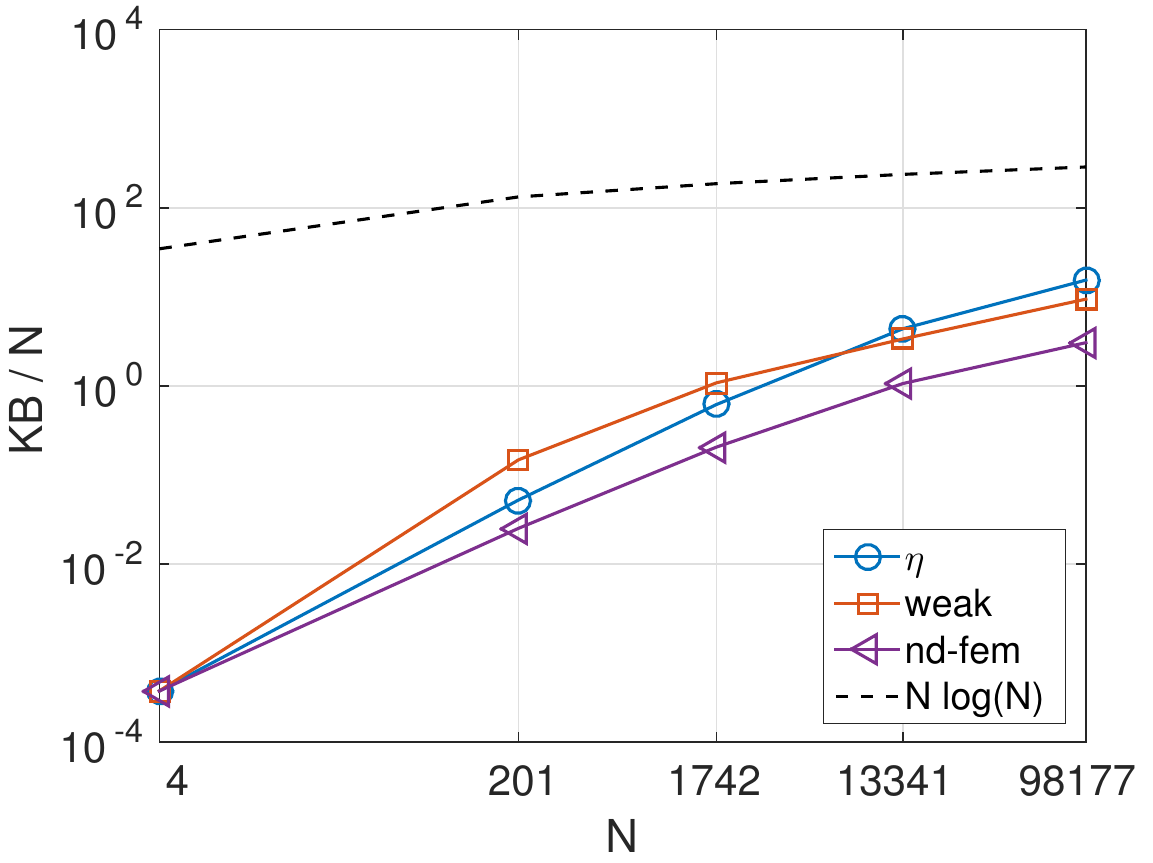}\quad
\includegraphics[width=0.35\textwidth]{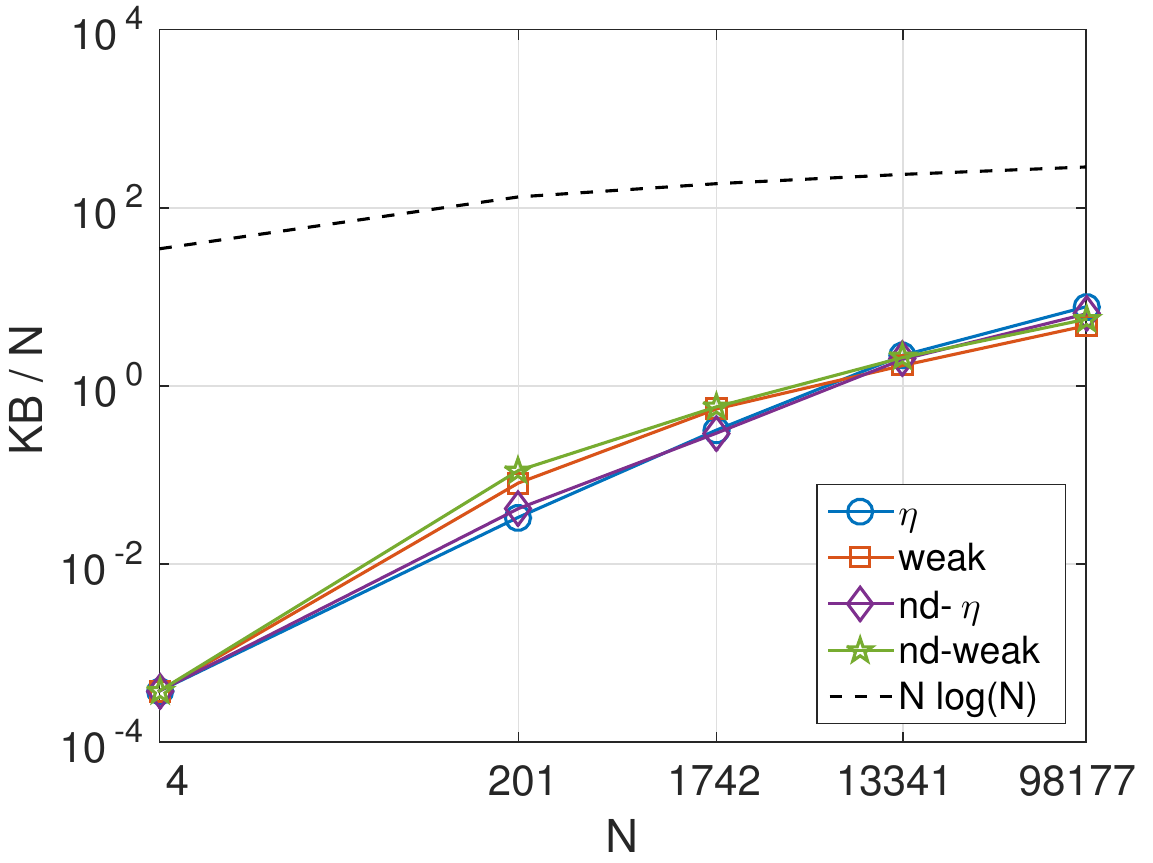}
\caption{\label{fig:bonestorage}Allocated storage per mesh degree of
freedom for different admissibility conditions. Storage for nonsymmetric
(left) and symmetric matrices (right). The allocated storage is
independent of the content of the matrix.}
\end{figure}

We are also interested in the quality of the approximate LU-decomposition
$\mathbf{A}\approx\hat{\mathbf{L}}\hat{\mathbf{U}}$.
We use a built-in function of the HLib to estimate the deviation 
$\hat{\mathbf{L}}\hat{\mathbf{U}}-\mathbf{A}$ in the spectral norm 
by ten power iterations, which is a good indicator for the approximation 
quality of the LU-decomposition of the finite element matrix. The estimated
errors are plotted in Figure~\ref{fig:deviation}. Note that the
observed behavior is in contrast to the behavior typically observed for
preconditioning, cf.,~e.g.,~\cite{Beb08}, since we do not increase the rank
with the number of unknowns. We can see that the 
LU-decomposition is most accurate in the all-$\eta$ case. Still, only
one iteration is needed in the iterative refinement in all cases. When it 
comes to computation times, Figure~\ref{fig:bonetimes} and 
Tables~\ref{tab:bonetrhs}, \ref{tab:bonetinv} and \ref{tab:bonetsolve} indicate that
all cases of admissibility under consideration might yield essentially 
linear complexity, although the asymptotic regime seems 
not to be reached in the considered levels of refinement.
Both, the weak admissibility condition and the nested 
dissection approach lead to considerable speed-ups, 
where the combination of these approaches, the 
nd-weak-case, seems to be the fastest approach.
Figure~\ref{fig:boneranks} illustrates the required 
average and maximal ranks needed for the computations, whereas
Figure~\ref{fig:bonestorage} illustrates the amount of storage needed 
per mesh degree of freedom. For reasons of performance, the HLib 
allocates the worst case scenario for the ranks. Thus, in the latter case,
only the different admissibilities for a single $\mathcal{H}$-matrix 
build from a binary and a nested dissection cluster tree have 
to be considered. In conclusion, the nested dissection clustering 
consumes less computation time and less storage for the LU-decomposition.

\begin{table}[hbt]
\centering
\begin{tabular}{|l|ccccc|}
\hline
Level       &          1 &          2 &          3 &          4 &          5   \\\hline\hline
all-$\eta$  &   0.000954 &   0.123520 &    8.07820 &    367.071 &   8,158.12	  \\
weak-FEM    &   0.001476 &   0.109641 &    8.12224 &    380.533 &   8,370.79   \\
all-weak    &   0.001113 &   0.113122 &    8.86575 &    434.289 &   9,464.19   \\
nd-$\eta$   &   0.001569 &   0.151018 &    7.54221 &    319.773 &   6,276.75   \\
nd-weak     &   0.000885 &   0.123694 &    8.07407 &    349.585 &   6,711.85   \\\hline
\end{tabular}
\caption{\label{tab:bonetrhs}Computation times in seconds to compute the data correlation $\mathcal{H}$-matrix $\mathbf{C}_f$ on the dumbbell geometry.}
\end{table}

\begin{table}[hbt]
\centering
\begin{tabular}{|l|ccccc|}
\hline
Level       &          1 &          2 &          3 &          4 &          5   \\\hline\hline
all-$\eta$  &   $2.2\cdot10^{-5}$ &   0.001274 &   0.339432 &    56.2521 &    2,806.5  \\
weak-FEM    &    $2.9\cdot10^{-5}$ &    0.00143 &   0.315344 &    17.7348 &    743.624  \\
all-weak    &    $2.2\cdot10^{-5}$ &   0.001831 &   0.316497 &    18.3358 &    746.364  \\
nd-$\eta$   &    $2.5\cdot10^{-5}$ &   0.000513 &   0.048170 &    4.17870 &    135.778  \\
nd-weak     &    $3.6\cdot10^{-5}$ &   0.000588 &   0.050896 &    4.05899 &    132.319  \\\hline
\end{tabular}
\caption{\label{tab:bonetinv}Computation times in seconds to compute the approximate
         LU-decomposition of the finite element matrix on the dumbbell geometry.}
\end{table}

\begin{table}[hbt]
\centering
\begin{tabular}{|l|ccccc|}
\hline
Level       &          1 &          2 &          3 &          4 &          5   \\\hline\hline
all-$\eta$  &    $5.2\cdot10^{-5}$ &   0.011115 &    4.68355 &   1,419.19 &   85,477.9	   \\
weak-FEM    &   0.000104 &   0.010592 &    2.64065 &    420.153 &   29,492.1   \\
all-weak    &    $5.4\cdot10^{-5}$ &   0.042098 &    14.5121 &    691.129 &   24,225.6   \\
nd-$\eta$   &   0.000102 &   0.039209 &    5.60769 &    443.310 &   21,390.0   \\
nd-weak     &    $4.4\cdot10^{-5}$ &   0.061530 &    7.67542 &    544.080 &   18,570.5   \\\hline
\end{tabular}
\caption{\label{tab:bonetsolve}Computation times in seconds for the iterative refinement 
on the dumbbell geometry.}
\end{table}

Having verified the convergence of our solver, we now want 
to consider different correlation lengths and different classes of smoothness.

\subsection{Small correlation lengths}
In the second part of the numerical experiments, we employ
correlation kernels with smaller correlation lengths and lower 
regularity such that low-rank approximations would become 
prohibitively expensive and sparse tensor product approaches would 
fail to resolve the concentrated measure.

\begin{figure}
\centering
\includegraphics[width=0.3\textwidth, clip=true, trim=180 540 180 100]{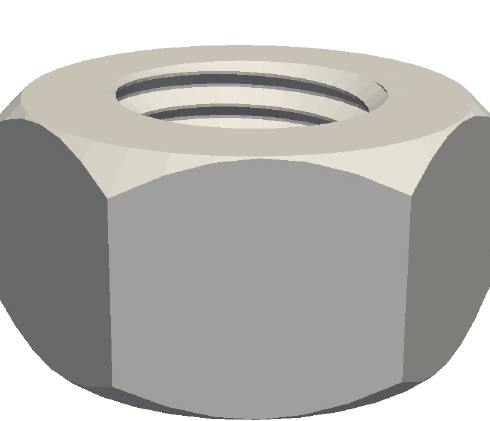}\qquad
\includegraphics[width=0.3\textwidth, clip=true, trim=180 540 180 100]{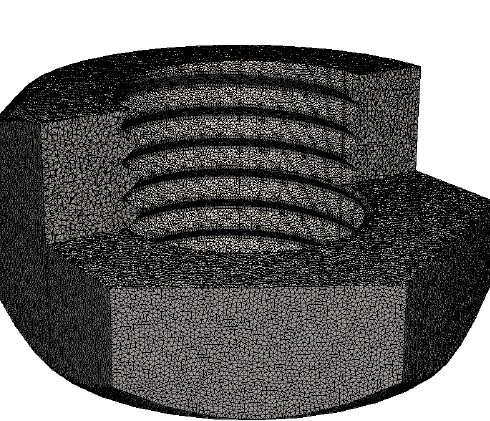}
\caption{\label{fig:nut}The screw-nut geometry (left) and its meshed cross section (right).}
\end{figure}

We consider the screw-nut geometry pictured in Figure~\ref{fig:nut} 
which is discretized by a mesh with 269,950 vertices, 197,480
mesh degrees of freedom, and a maximal element diameter of 
$h/\diam(D)\approx 0.0225$, yielding a matrix equation with
$3.90\cdot10^{10}$ unknowns. We choose $\mathcal{L} = -\Delta$ 
in \eqref{eq:coreq} and either the Gaussian kernel as input correlation 
$\Cor _f$, i.e.,
\[
\Cor _f(\mathbf{x},\mathbf{y}) = \frac{1}{\ell}\exp\bigg(-\frac{\|\mathbf{x}-\mathbf{y}\|^2}{2\ell^2}\bigg),
\]
or the exponential kernel, i.e.,
\[
\Cor _f(\mathbf{x},\mathbf{y}) = \frac{1}{\ell}\exp\bigg(-\frac{\|\mathbf{x}-\mathbf{y}\|}{\ell}\bigg).
\]
Herein, $\ell>0$ denotes the correlation length.

In the following, we want to demonstrate that the presented method 
is well suited for small correlation lengths $\ell$. We therefore 
choose the correlation lengths
\[
\ell \in\bigg\{\frac{\diam (D)}{1},\frac{\diam(D)}{2},\frac{\diam(D)}{4},
\frac{\diam(D)}{8},\frac{\diam(D)}{16},\frac{\diam(D)}{32}\bigg\}
\]
for both, the Gaussian kernel and the exponential kernel, and compute
the corresponding correlation of the solution $\Cor _u$ of \eqref{eq:coreq}.

In our first experiment, we use the nd-weak case, as the 
previous section has shown that it is more memory efficient and
has superior computation times.
The computation time for the assembly of the prescribed correlation is
around 20,000 seconds and the computation time of the approximate 
LU-decomposition is around 400 seconds,
whereas the computation times for the iterative refinement are contained 
in Table~\ref{tab:nutsolve}.

\begin{table}
\centering
\scalebox{0.88}{
\begin{tabular}{|c|c|cccccc|}
\hline
                             & $\ell/\diam(D)$ & 1 &      1/2 &      1/4 &      1/8 &     1/16 &     1/32 \\\hline\hline
\multirow{2}{*}{Exponential} & nd-weak  & 51,656.7 & 53,011.0 & 52,876.5 & 51,459.2 & 49,838.2 & 51,524.6\\
                             & all-weak & 77,784.5 & 79,101.8 & 79,155.3 & 79,155.3 & 76,952.6 & 72,256.9\\\hline\hline
\multirow{2}{*}{Gaussian}    & nd-weak  & 47,921.8 & 50,644.0 & 50,819.5 & 51,753.7 &      --- &      --- \\
                             & all-weak & 73,405.4 & 74,877.0 & 75,165.7 & 68,222.8 & 72,259.4 & 75,070.4 \\\hline
\end{tabular}
}
\caption{\label{tab:nutsolve}Computation times in seconds 
for the nd-weak case and for the all-weak case
for the iterative refinement on the screw-nut geometry for 
the exponential kernel and the Gaussian kernel with 
different correlation lengths.}
\end{table}

We do not tabulate the computation times for the Gaussian kernel for the
correlation lengths $\ell/\diam(D)=1/16$ and $\ell/\diam(D)=1/32$
since the iterative refinement does not converge to the prescribed 
tolerance. In all other cases, the iterative refinement needs only 
one iteration.

Repeating the computations in the two problematic cases
with increased $k_\eta$ or in the nd-$\eta$ instead
of the nd-weak case does also not lead to convergence. However, repeating all
computations in the all-weak case resolves the issue, as the computation
times in Table~\ref{tab:nutsolve} show. In the all-weak case, the
computation time for the prescribed correlation is again around 20,000 seconds
and the computation time for the approximate LU-decomposition is around 1,700
seconds. The iterative refinement needs again one iteration in all tabulated
cases.

\begin{figure}
\centering
\begin{minipage}{0.28\textwidth}
\centering
\includegraphics[width=\textwidth, clip=true, trim=160 520 30 100]{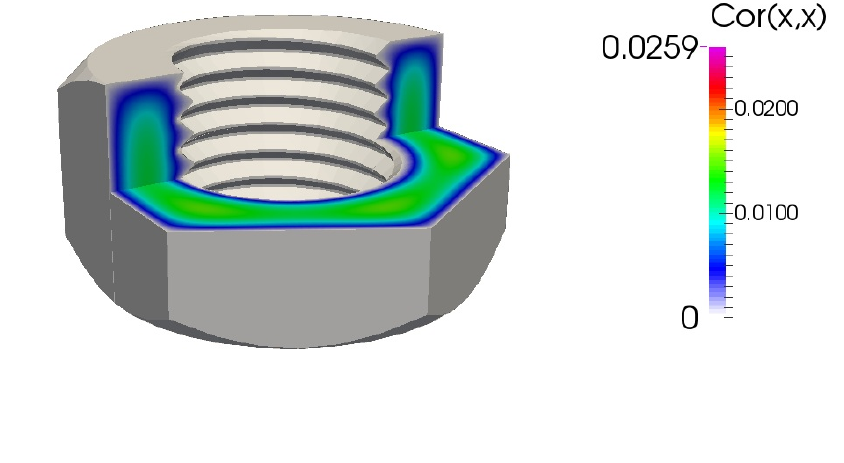}\\
$\ell =\frac{\diam (D)}{1}$
\end{minipage}\quad
\begin{minipage}{0.28\textwidth}
\centering
\includegraphics[width=\textwidth, clip=true, trim=160 520 30 100]{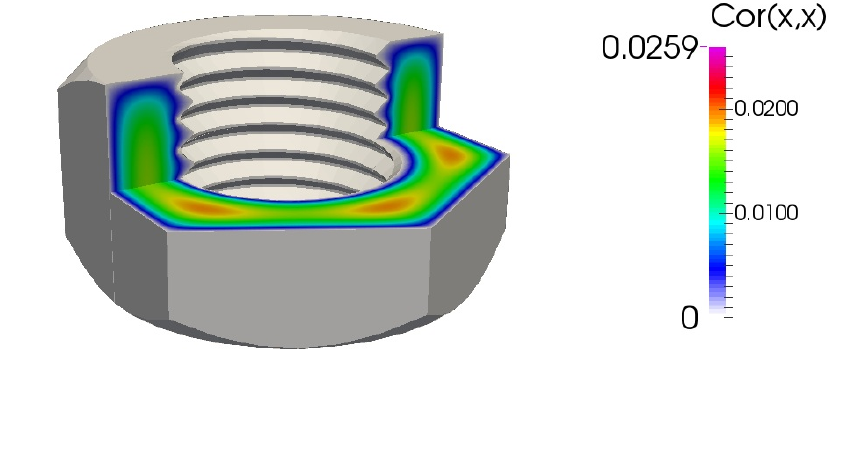}\\
$\ell =\frac{\diam (D)}{2}$
\end{minipage}\quad
\begin{minipage}{0.28\textwidth}
\centering
\includegraphics[width=\textwidth, clip=true, trim=160 520 30 100]{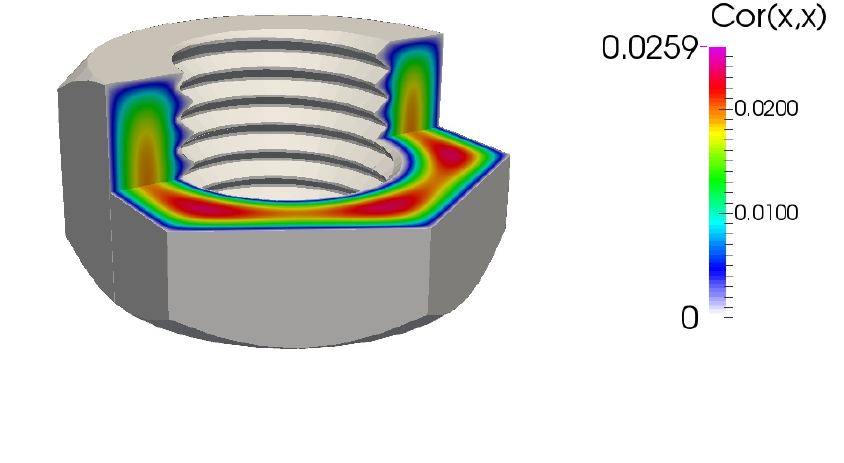}\\
$\ell =\frac{\diam (D)}{4}$
\end{minipage}\\
\medskip
\begin{minipage}{0.28\textwidth}
\centering
\includegraphics[width=\textwidth, clip=true, trim=160 520 30 100]{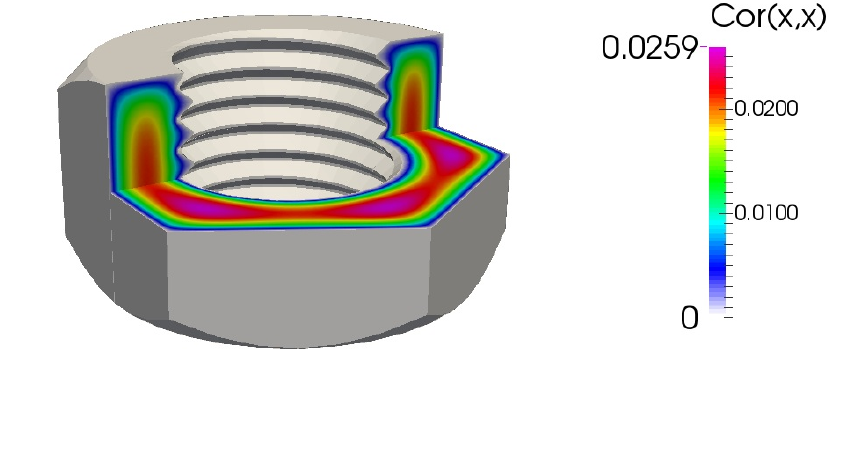}\\
$\ell =\frac{\diam (D)}{8}$
\end{minipage}\quad
\begin{minipage}{0.28\textwidth}
\centering
\includegraphics[width=\textwidth, clip=true, trim=160 520 30 100]{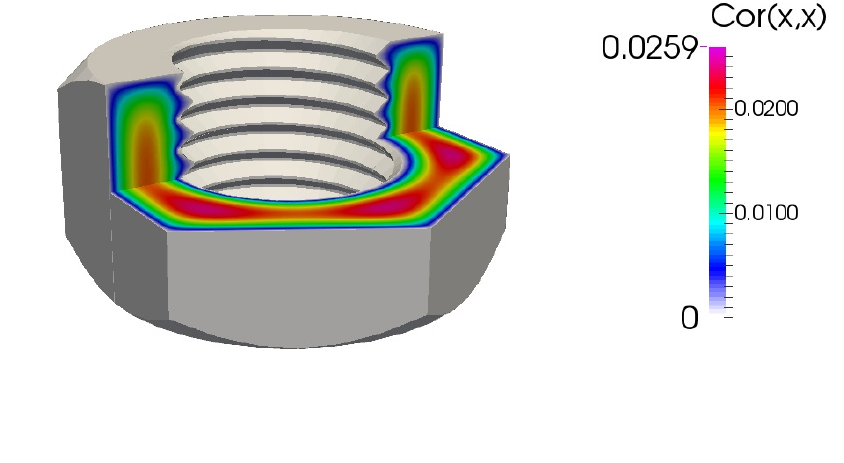}\\
$\ell =\frac{\diam (D)}{16}$
\end{minipage}\quad
\begin{minipage}{0.28\textwidth}
\centering
\includegraphics[width=\textwidth, clip=true, trim=160 520 30 100]{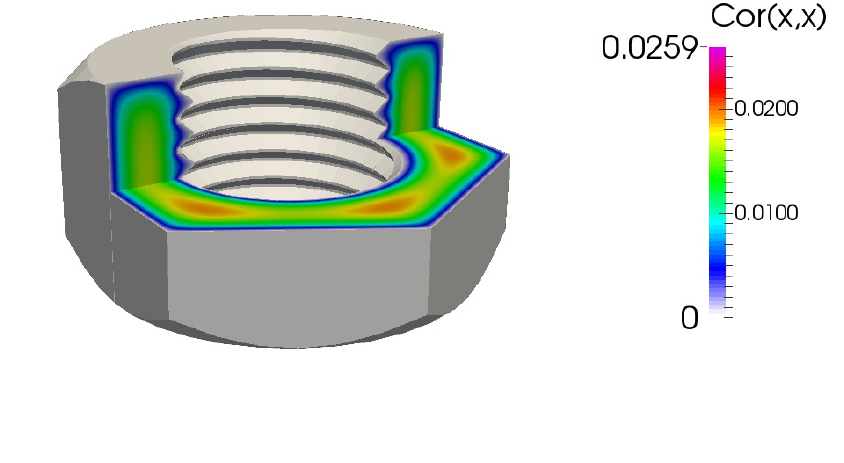}\\
$\ell =\frac{\diam (D)}{32}$
\end{minipage}\\
\caption{\label{fig:nutexpcross}Cross sections of the correlation of the solution through the screw-nut 
geometry for the exponential kernel with different correlation lengths $\ell$.}
\end{figure}

\begin{figure}
\centering
\begin{minipage}{0.28\textwidth}
\centering
\includegraphics[width=\textwidth, clip=true, trim=160 520 30 100]{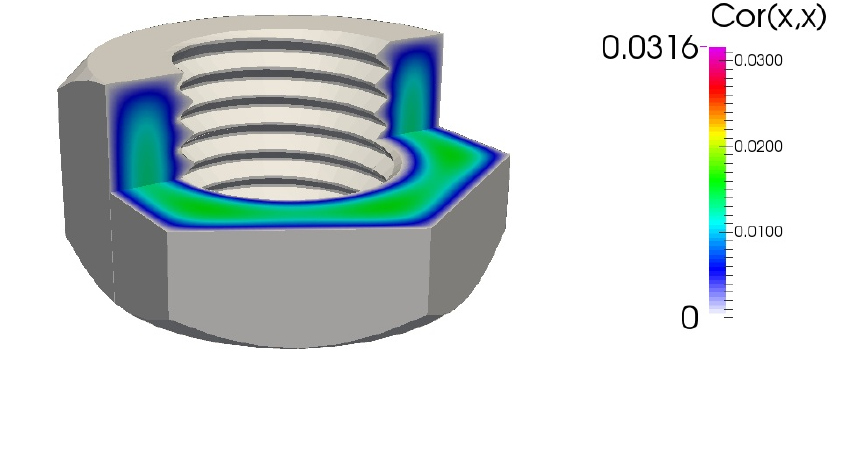}\\
$\ell =\frac{\diam (D)}{1}$
\end{minipage}\quad
\begin{minipage}{0.28\textwidth}
\centering
\includegraphics[width=\textwidth, clip=true, trim=160 520 30 100]{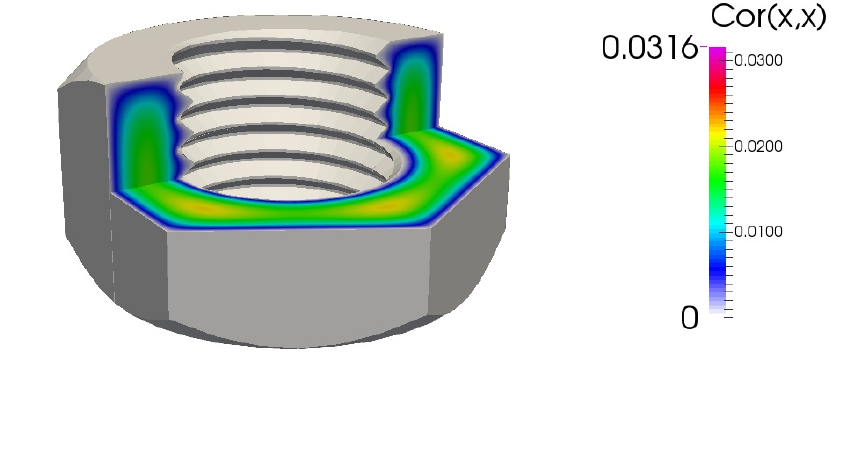}\\
$\ell =\frac{\diam (D)}{2}$
\end{minipage}\quad
\begin{minipage}{0.28\textwidth}
\centering
\includegraphics[width=\textwidth, clip=true, trim=160 520 30 100]{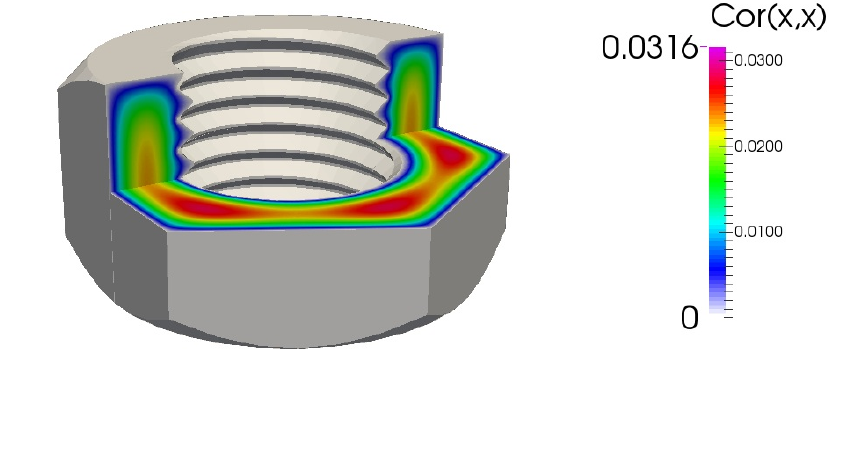}\\
$\ell =\frac{\diam (D)}{4}$
\end{minipage}\\
\medskip
\begin{minipage}{0.28\textwidth}
\centering
\includegraphics[width=\textwidth, clip=true, trim=160 520 30 100]{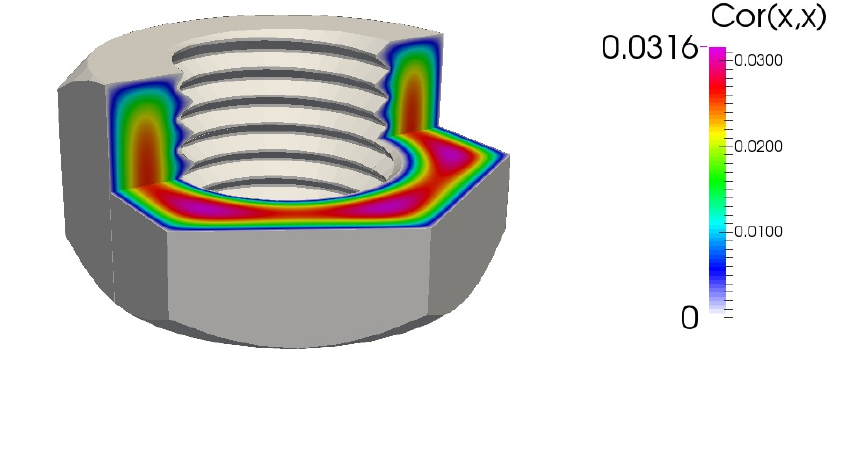}\\
$\ell =\frac{\diam (D)}{8}$
\end{minipage}\quad
\begin{minipage}{0.28\textwidth}
\centering
\includegraphics[width=\textwidth, clip=true, trim=160 520 30 100]{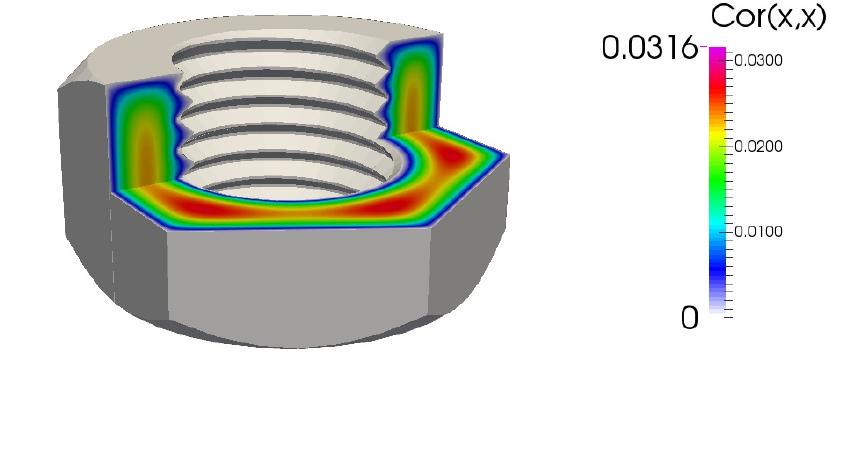}\\
$\ell =\frac{\diam (D)}{16}$
\end{minipage}\quad
\begin{minipage}{0.28\textwidth}
\centering
\includegraphics[width=\textwidth, clip=true, trim=160 520 30 100]{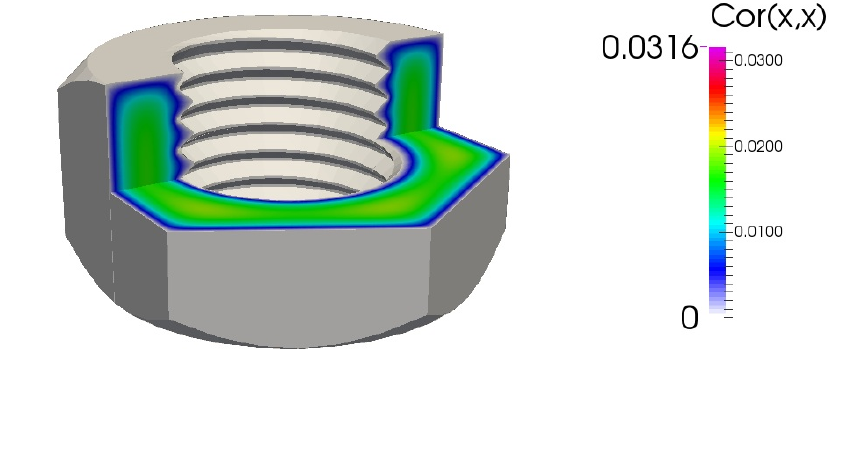}\\
$\ell =\frac{\diam (D)}{32}$
\end{minipage}\\
\caption{\label{fig:nutgausscross}Cross sections of the correlation of the solution through the screw-nut geometry 
for the Gaussian kernel with different correlation lengths $\ell$.}
\end{figure}

The cross sections found in Figure~\ref{fig:nutexpcross} 
illustrate the different behaviour of the correlation's trace 
$\Cor _u|_{\mathbf{x}=\mathbf{y}}$ for the different correlation
lengths in case of the exponential kernel. The related results
for the Gaussian kernel are presented in Figure~\ref{fig:nutgausscross}.
It seems that there occours a mass defect in the correlation lengths
$\ell/\diam(D)=1/16$ and $\ell/\diam(D)=1/32$. This could be due to the fact
that the mesh size of the finite element method is not able to resolve the
correlation length properly. Nevertheless, the computation times are
independent of $\ell$, even if the underlying finite element method cannot
resolve the correlation length. Moreover, the nested dissection clustering
technique can lead to a speed-up, while the binary clustering technique
seems to be more robust.

\section{Conclusion}\label{sec:conclusion}
We considered the solution of strongly elliptic partial differential 
equations with random load by means of the finite element method.
Approximating the full tensor approach by means of $\mathcal{H}$-matrices,
we employed the $\mathcal{H}$-matrix technique to efficiently discretize 
the non-local correlation kernel of the data and approximate the LU-decomposition
of the finite element stiffness matrix. The corresponding $\mathcal{H}$-matrix equation 
has then been efficiently solved in essentially linear complexity by the 
$\mathcal{H}$-matrix arithmetic. 

Compared to sparse tensor product or low-rank approximations, the 
proposed method does not suffer in case of shortly correlated 
data from large constants in the complexity estimates or the lack 
of resolution of the roughness. This has been shown by numerical 
experiments on a non-trivial three-dimensional geometry. Indeed,
neither the  computation times nor the storage requirements do increase
for  correlation kernels with short correlation length. It was moreover 
demonstrated that the use of the weak admissibility condition
for the partition of the $\mathcal{H}$-matrix improves the constants in
the computational complexity without having a significant impact to the
solution accuracy. The use of a nested dissection clustering strategy can
additionally lead to a speed-up of the computations and save storage,
whereas the binary clustering strategy seems to be the more robust approach.

\section*{Acknowledgement}
The authors would like to thank Steffen B\"orm and the anonymous referees 
for valuable and helpful remarks on the weak admissibility and the
nested dissection $\mathcal{H}$-LU decomposition.

\bibliographystyle{amsplain}
\bibliography{bibl}
\end{document}